# A generalized stochastic control problem of bounded noise process under ambiguity arising in biological management

Hidekazu Yoshioka[1, 2, *] and Motoh Tsujimura[3]

[1] Assistant Professor, Graduate School of Natural Science and Technology, Shimane University, Nishikawatsu-cho 1060, Matsue, 690-8504 Japan

[2] Center member, Fisheries Ecosystem Project Center, Shimane University, Nishikawatsu-cho 1060, Matsue, 690-8504, Japan

[3] Professor, Faculty of Commerce, Doshisha University, Karasuma-Higashi-iru, Imadegawa-dori, Kamigyo-ku, Kyoto, 602-8580 Japan

* Corresponding author

E-mail: yoshih@life.shimane-u.ac.jp

**Abstract**

The objectives and contributions of this paper are mathematical and numerical analyses of a stochastic control problem of bounded population dynamics under ambiguity, an important but not well-studied problem, focusing on the optimality equation as a nonlinear degenerate parabolic partial integro-differential equation (PIDE). The ambiguity comes from lack of knowledge on the continuous and jump noises in the dynamics, and its optimization appears as nonlinear and nonlocal terms in the PIDE. Assuming a strong dynamic programming principle for continuous value functions, we characterize its solutions from both viscosity and distribution viewpoints. Numerical computation focusing on an ergodic case are presented as well to complement the mathematical analysis.

**AMS Classification**

93E20: Optimal stochastic control

35K65: Degenerate parabolic equations

49L25: Viscosity solutions

35D30: Weak solutions

65M06: Finite difference methods

## 1. Introduction

Many biological processes are stochastic and bounded because of some nonlinearity as found in the classical logistic dynamics [1]. Diffusion processes having bounded ranges, the simplest one being the stochastic logistic-type model [2], are commonly found as efficient stochastic models of biological processes. Jump-diffusion processes and regime-switching diffusion processes have been employed as advanced alternatives to the models based on bounded diffusion processes. Epidemiological models under complete information [3, 4] and incomplete information [5] have been theoretically analyzed with bounded diffusion processes. Yoshioka [6] discussed an infinite-horizon problem of algae population management in rivers based on controlled jump processes. A multi-dimensional sediment-algae interaction problem has also been considered Yoshioka et al. [7]. The model of Lungu and Øksendal [2] has been used as a growth curve of fish [8]. Often the bounded noises are more adequate than the unbounded ones in physics and engineering as well [9]. Furthermore, some stochastic epidemic models assume bounded population dynamics [10].

From a management viewpoint, modeling a target phenomenon is often not the goal, but is a path toward its optimization. Stochastic optimal control [11] provides a platform for modeling, analysis, and optimization of a variety of system dynamics in the real world. A conventional deterministic optimal control problem can be seen as a non-stochastic counterpart of a stochastic control problem. The application examples are diverse, ranging from biology [12] to mathematical physics [13], environmental sciences [14], and social sciences [15] even if only to recent research. Stochastic models in applications often involve some ambiguities, or equivalently misspecifications, due to our lack of knowledge of the target system dynamics.

The multiplier robust control [16] is the mathematical concept for efficiently handling this issue based on a differential game between the decision-maker and the opponent player, called nature, representing worst-case ambiguity. The modeling approach based on the multiplier robust control has been successfully applied to problems with unbounded and often analytically tractable dynamics in finance [17], energy and economics [18, 19], and insurance [20, 21]. Application of the approach to environmental problems has been carried out as well, demonstrating its versatility [22]. However, application of the multiplier robust control to problems with bounded diffusion and jump-diffusion processes are much rarer despite their relevance in applied problems like biological management, except for several recent research [6, 23, 24]. Especially, to the authors' best of knowledge, mathematical analysis like solvability and regularity of the optimality equation in such a case has not been carried out so far. This is the motivation of this paper.

The objective of this paper is to carry out mathematical analysis of the optimality equation, which is a Hamilton-Jacobi-Bellman-Isaacs (HJBI) equation associated with a model bounded jump-diffusion process. Using HJBI equations in stochastic control has a clear advantage that is the optimal controls are found without resorting to statistical simulation based on Monte-Carlo type methods. Our model problem considers single-species population dynamics as a generalization of the stochastic logistic model [2], and is driven by ambiguous continuous and jump noises. The HJBI equation is a degenerate parabolic partial integro-differential equation (PIDE) having nonlinear and nonlocal terms. The model generalizes several population dynamics models, and the HJBI has not been found in the literature. We demonstrate that the HJBI equation can be solvable in both continuous viscosity solutions [25, 26] and distribution [27, 28] senses. The latter seems to be less common for analyzing HJBI equations. Note that we avoid using the conventional term "weak solution" for the distribution solution since a viscosity solution is also a weak solution of the HJBI equation as well in a different sense.

Both the viscosity and distribution approaches have pros and cons. An advantage of the viscosity approach is its ability to analyze existence and uniqueness of continuous, but possibly non-smooth solutions, in a unified manner under relatively weak regularity conditions of the coefficients. Often even the continuity assumption of solutions is unnecessary. On the other hand, its disadvantage is that analyzing regularity of the solutions, such as their differentiability, is not always easy. In contrast, the distribution approach requires relatively stronger regularity conditions on the coefficients, but solutions obtained belong to some weighted Sobolev spaces that can guarantee some boundedness of the solutions and their partial derivatives. Linkages between the viscosity solutions and the distribution solutions are nontrivial even for linear problems without nonlocal terms [29]. In this paper, we employ both approaches to obtain deeper understanding of solutions to the HJBI equation. Assuming a strong dynamic programming principle as in the past research [16,17, 20-22], we demonstrate that the HJBI equation has at least the two characterizations based on the separate approaches. Furthermore, we show that a simple finite difference scheme verifying the desirable mathematical properties [30] can numerically handle the HJBI equation in a stable, monotone, and consistent manner. Demonstrative computational results are presented as well.

The rest of this paper is organized as follows. Section 2 presents the model problem and derives the HJBI equation. Section 3 is devoted to unique solvability and regularity results of the HJBI equation. A part of Section 3 focuses on its numerical computation. Section 4 concludes this paper and presents our future perspectives.

## 2. Mathematical model

**2.1 Mathematical setting**

We consider population dynamics of single species, which is some animal population or some index like water quality, in a habitat or a region [6, 23, 24, 31]. All the stochastic processes appearing this paper are defined in a standard complete probability space as in the conventional cases [11]. The decision-maker, the observer, tries to control the dynamics so that a performance index is minimized. On the other hand, nature as an opponent maximizes the performance index.

The time is denoted as $t \geq 0$ and the population is set as a càdlàg variable $X_t$ at time $t$ with the initial condition $x \in \Omega = [0,1]$. The 1-D standard Brownian motion is denoted as $B_t$ at time $t$. The standard compound Poisson process with the intensity $\nu_i > 0$ and the jump probability density $g_i$ is denoted as $N_t^{(i)}$ at time $t$ ($i = 1,2$). The processes $(B_t)_{t\geq 0}$, $\left(N_t^{(1)}\right)_{t\geq 0}$, and $\left(N_t^{(2)}\right)_{t\geq 0}$ are assumed to be mutually independent. For each $i = 1,2$, the jump magnitude $z_i$ of $N_t^{(i)}$ is assumed to be valued in a compact set in the open set $(0,1)$. Therefore, each $g_i$ is compactly support in $(0,1)$ and satisfies $\int_0^1 g_i(z)\mathrm{d}z = 1$ and $0 < \int_0^1 \frac{g_i(z)}{1-z}\mathrm{d}z < +\infty$. A natural filtration generated by $(B_t)_{t\geq 0}$, $\left(N_t^{(1)}\right)_{t\geq 0}$, and $\left(N_t^{(2)}\right)_{t\geq 0}$ is denoted as $\mathcal{F} = (\mathcal{F}_t)_{t>0}$.

**2.2 Nominal model**

Our SDE is a minimum logistic-type SDE considering growth, migration/immigration, and continuous/jump noises of controlled single-species population. The model SDE without ambiguity is formulated in the Itô's sense as

$$\begin{aligned}\mathrm{d}X_t = {} & a(X_t)\left(r(q_t)\mathrm{d}t + \sigma \mathrm{d}B_t\right) - \gamma_0 X_t \mathrm{d}t + \gamma_1 (1 - X_t)\mathrm{d}t \\ & - X_{t-0}\mathrm{d}N_t^{(1)} + (1 - X_{t-0})\mathrm{d}N_t^{(2)}\end{aligned}, \quad t > 0, \tag{1}$$

which is a generalization of the logistic model [2, 4]. In the right-hand side of (1), the first term is the continuous stochastic growth term with the density-dependent growth rate $a : \Omega \to \mathbb{R}$, where $r > 0$ and $\sigma > 0$ are the magnitudes of the deterministic and stochastic parts of the growth, respectively. The deterministic growth rate $r$ is controlled through the intervention $(q_t)_{t>0}$ adapted to $\mathcal{F}$ having the compact range $Q = [0, q_{\max}]$ with a constant $q_{\max} > 0$. We assume $r$ is Lipschitz continuous in $Q$. Set $r_{\max} = \max_{q\in Q} r(q) > 0$. The second term represents the decaying and/or migration with the intensity $\gamma_0 \geq 0$ and the third term represents the immigration with the intensity $\gamma_1 \geq 0$. The fourth and fifth terms are the jump decrease and increase of the population, respectively. The SDE (1) is a simple generalized mathematical model of stochastically varying bounded population dynamics having a martingale term (the term with $\mathrm{d}B_t$), a decreasing term (the term with $\mathrm{d}N_t^{(1)}$), and an increasing term (the term with $\mathrm{d}N_t^{(2)}$). We employ **Assumption 1** for well-posedness of the SDE (1) following the logistic type framework of stochastic population dynamics [2, 32].

***Assumption 1***

$a(0)=a(1)=0$, $a>0$ in $(0,1)$, and $a$ is Lipschitz continuous in $\Omega$.

***Remark 1***

By **Assumption 1**, the SDE (1) is uniquely solvable in the path-wise strong sense, and the solution $(X_t)_{t\geq 0}$ is a.s. in $\Omega$. This statement follows from a contradiction argument [2] combined with the Itô's formula for jump-diffusion processes.

***Remark 2***

The fourth and fifth terms of (1) are chosen based on a formal symmetry to simplify the model. In fact, they satisfy $\mathrm{d}X_t = -z_1 X_{t-0}$ and $\mathrm{d}(1-X_t) = -z_2(1-X_{t-0})$ at jumps, respectively.

### 2.3 Ambiguous model

The SDE (1) is extended to an ambiguous counterpart based on the multiplier robust formalism [16]. The model ambiguity in this framework is represented as distortions of the reference probability measure denoted as $\mathbb{P}$. Formally, we consider the distorted (misspecified) processes generated from $(B_t)_{t>0}$, $\left(N_t^{(1)}\right)_{t>0}$, and $\left(N_t^{(2)}\right)_{t>0}$: $\hat{B}_t = -\int_0^t \lambda_s \mathrm{d}s + B_t$, and $\hat{N}_t^{(i)}$ is the compound Poisson processes having the jump probability density $g_i$ and the intensity $v_i\theta_t^{(i)}$ at time $t$ ($i=1,2$). $(\lambda_t)_{t>0}$ and $\left(\theta_t^{(i)}\right)_{t>0}$ ($i=1,2$) are adapted to $\mathcal{F}$ and valued in $\Lambda = [-\lambda_{\max}, \lambda_{\max}]$ and $\Theta = [0, \theta_{\max}]$ with sufficiently large constants $\lambda_{\max}, \theta_{\max} > 0$. We accordingly obtain a non-ambiguous model when $\theta_t^{(i)} \equiv 1$ and $\lambda_t \equiv 0$. We consider $0 \ln 0 = 0$.

Assume $\int_0^t \left(\theta_s^{(i)} \ln \theta_s^{(i)} + 1 - \theta_s^{(i)}\right)\mathrm{d}s < +\infty$ a.s. and $\exp\left(\int_0^t \left(\theta_s^{(i)} \ln \theta_s^{(i)} + 1 - \theta_s^{(i)}\right)\mathrm{d}s\right) < +\infty$ in the mean for $t>0$ and $i=1,2$. Then, set a new probability measure $\mathbb{Q}$, which is related to the original probability measure $\mathbb{P}$ through the Radon-Nikodym derivative (See, Hansen and

Sargent [16] and Hu and Wang [33])

$$\frac{\mathrm{d}\mathbb{Q}}{\mathrm{d}\mathbb{P}} = \Xi_t^{(0)}\Xi_t^{(1)}\Xi_t^{(2)}, \tag{2}$$

where

$$\Xi_t^{(0)} = \exp\left(-\frac{1}{2}\int_0^t \lambda_s^2 \mathrm{d}s + \int_0^t \lambda_s \mathrm{d}B_s\right), \tag{3}$$

$$\Xi_t^{(i)} = \exp\left(\int_0^t \ln\theta_s^{(i)} \mathrm{d}N_s^{(i)} + \int_0^t\int_0^1 \left(1-\theta_s^{(i)}\right)\nu_i g_i(z)\mathrm{d}z\mathrm{d}s\right) \quad (i=1,2). \tag{4}$$

Then, $\mathbb{Q}$ is considered as a probability measure under distortions due to the noise ambiguity. Now, the ambiguous counterpart of the SDE (1) under $\mathbb{Q}$ is formulated as

$$\begin{aligned}\mathrm{d}X_t &= a(X_t)\left(r(q_t)\mathrm{d}t + \sigma\mathrm{d}\hat{B}_t\right) - \gamma_0 X_t \mathrm{d}t + \gamma_1(1-X_t)\mathrm{d}t \\ &\quad + \sigma\lambda_t a(X_t)\mathrm{d}t - X_{t-0}\mathrm{d}\hat{N}_t^{(1)} + (1-X_{t-0})\mathrm{d}\hat{N}_t^{(2)}\end{aligned}, \quad t>0. \tag{5}$$

There exists an additional distortion term in the SDE (5). All the expectations appearing below are defined in the sense of the probability measure $\mathbb{Q}$.

***Remark 3***

By **Assumption 1**, the SDE (5) is also uniquely solvable in the path-wise strong sense, and the solution $X_t$ is almost surely bounded in the compact domain $\Omega$ for $t>0$.

### 2.4 HJBI equation

The performance index is an expectation to be minimized by the decision-maker while maximized by nature:

$$\phi\left(t,x;q,\lambda,\theta^{(1)},\theta^{(2)}\right) = \mathbb{E}\left[\int_t^T \left(f(X_s) + h(q_s) - g_0(\lambda_s) - \sum_{i=1}^2 g_i\left(\theta_s^{(i)}\right)\right)\mathrm{d}s\right], \tag{6}$$

where $X_{t-0} = x$. $f$ is a non-negative Lipschitz continuous function in $D$ with $f_{\max} = \max_{\Omega} f > 0$ representing the net disutility caused by the population, $h$ is a non-negative and Lipschitz continuous function in $Q$ as a unit-time control cost. Each $g^{(i)}$ is the entropic penalization term, which is defined following the multiplier robust approach [16, 19, 33] as

$$g_0(\lambda_s) = \frac{\lambda_s^2}{2\psi_0}, \quad g_i\left(\theta_s^{(i)}\right) = \frac{\nu_i}{\psi_i}\left(\theta_s^{(i)}\ln\theta_s^{(i)} + 1 - \theta_s^{(i)}\right) \quad (i=1,2), \tag{7}$$

where $\psi_0$, $\psi_1$, $\psi_2$ are the positive constants serving as the ambiguity-aversion parameters of the ambiguity on the processes $(B_t)_{t\geq 0}$, $\left(N_t^{(1)}\right)_{t\geq 0}$, $\left(N_t^{(2)}\right)_{t\geq 0}$, respectively. Each coefficient in

(7) is proportional to the relative entropy between the true and distorted models, and thus penalizes deviations between them. The decision-maker sets $\psi_i$, and its larger value means more ambiguity-aversion to the corresponding noise. The case $\psi_i \to +0$ means the ambiguity-neutral case that the decision-maker ignores the ambiguity.

The admissible set $\mathcal{A}$ of the process $(q_t)_{t>0}$ as a collection of the control variables chosen by the decision-maker is set as

$$\mathcal{A} = \left\{ (q_t)_{t>0} \left| \begin{array}{l} q_t \text{ is measurable, adapted to } \mathcal{F}_t, \\ \text{and valued in } Q \text{ for } t>0. \end{array} \right. \right\}. \tag{8}$$

The admissible set $\mathcal{B}$ of the triplet $\left(\lambda_t, \theta_t^{(1)}, \theta_t^{(2)}\right)_{t>0}$ as a collection of the control variables chosen by nature, the opponent, is set as

$$\mathcal{B} = \left\{ \left(\lambda_t, \theta_t^{(1)}, \theta_t^{(2)}\right)_{t>0} \left| \begin{array}{l} \left(\lambda_t, \theta_t^{(1)}, \theta_t^{(2)}\right) \text{ is measurable, adapted to } \mathcal{F}_t, \\ \text{and valued in } \Lambda \times \Theta \times \Theta \text{ for } t>0. \end{array} \right. \right\}. \tag{9}$$

The value function is set as the worst-case minimum value of the performance index

$$\Phi(t,x) = \inf_{q \in \mathcal{A}} \sup_{\lambda, \theta^{(1)}, \theta^{(2)} \in \mathcal{B}} \phi\left(t, x; \lambda, q, \theta^{(1)}, \theta^{(2)}\right), \tag{10}$$

which can be seen as an upper-value of a zero-sum game of a worst-case minimization type [16]. Clearly, we have the upper- and lower-bounds

$$0 \le \Phi(t,x) \le (T-t) f_{\max}, \tag{11}$$

suggesting that the limit $\lim_{T \to +\infty} T^{-1} \Phi(t,x)$, namely the ergodic limit, exists for each $(t,x)$ and $0 \le T^{-1}\Phi(t,x) \le f_{\max}$. This point is numerically explored in Section 3.

The optimizer of the right-hand side of (10) in the framework of Markov control [11] is called the optimal controls, and are denoted as $q_t^*, \lambda_t^*, \theta_t^{(1)*}, \theta_t^{(2)*}$, $t>0$. The goal of our problem is to find the optimal controls.

The HJBI equation is the governing PIDE of $\Phi$. Set $H: \Omega \times \mathbb{R} \to \mathbb{R}$ as

$$H(x,p) = f(x) - \max_{q \in Q}\left\{ -r(q) a(x) p - h(q) \right\} - \min_{\lambda \in \Lambda}\left\{ -\sigma \lambda a(x) p + \frac{\lambda^2}{2\psi_0} \right\}, \quad x \in \Omega, \quad p \in \mathbb{R}. \tag{12}$$

In addition, set $K, M$ acting on generic functions in $\Omega$ as

$$K\Phi(x) = \min_{\theta \in \Theta}\left\{ v_1 \theta \left( \Phi(x) - \int_0^1 g_0(z) \Phi((1-z)x) \mathrm{d}z \right) + \frac{v_1}{\psi_1}(\theta \ln \theta + 1 - \theta) \right\}, \quad x \in \Omega, \tag{13}$$

$$M\Phi(x)=\min_{\theta\in\Theta}\left\{\nu_2\theta\left(\Phi(x)-\int_0^1 g_0(z)\Phi(z+(1-z)x)\mathrm{d}z\right)+\frac{\nu_2}{\psi_2}(\theta\ln\theta+1-\theta)\right\},\quad x\in\Omega. \quad (14)$$

Their right-hand sides are further calculated as (See, **Remark 6**)

$$K(\Phi)(x)=\frac{\nu_1}{\psi_1}\left(1-e^{-\psi_1\left(\Phi(x)-\int_0^1 g_1(z)\Phi((1-z)x)\mathrm{d}z\right)}\right), \quad (15)$$

$$M(\Phi)(x)=\frac{\nu_2}{\psi_2}\left(1-e^{-\psi_2\left(\Phi(x)-\int_0^1 g_2(z)\Phi(z+(1-z)x)\mathrm{d}z\right)}\right), \quad (16)$$

respectively.

Set $D=[0,T)\times\Omega$. Assuming a strong dynamic programming principle for continuous value functions, our HJBI equation is

$$-\frac{\partial\Phi}{\partial t}-\frac{\sigma^2}{2}a(x)^2\frac{\partial^2\Phi}{\partial x^2}-\left[-\gamma_0 x+\gamma_1(1-x)\right]\frac{\partial\Phi}{\partial x}=H\left(x,\frac{\partial\Phi}{\partial x}\right)-K\Phi-M\Phi \quad \text{in } D \quad (17)$$

with the terminal condition $\Phi(T,x)=0$ for $x\in\Omega$. The operators $K,M$ in (17) act only on the second arguments. The right-hand sides are nonlinear as well as nonlocal. The Isaacs condition holds (17) since the maximization and minimization in $H$ can be carried out separately.

***Remark 4***

Continuity of the value function $\Phi$ on $D$ follows from the uniform continuity of the performance index $\phi$ on $D$ independent of the choice of control variables. This follows from the Lipchitz continuity of the coefficients of the SDEs and the Lipschitz continuity of $f$.

***Remark 5***

One may consider time-dependent $a$, $r$, and $\sigma$. The mathematical analysis results in Section 3 holds true in this case as well if their dependence on $r$ is smooth.

***Remark 6***

In the framework of Markov control, the optimal controls as functions in $D$ are found as

$$q^*(t,x)=\arg\max_{q\in Q}\left\{-r(q)a(x)\frac{\partial\Phi}{\partial x}-h(q)\right\}, \quad (18)$$

$$\lambda^*(t,x)=\arg\min_{\lambda\in\Lambda}\left\{-\sigma\lambda a(x)\frac{\partial\Phi}{\partial x}+\frac{\lambda^2}{2\psi_0}\right\}, \quad (19)$$

$$\theta^{(1)*}(t,x) = e^{-\psi_1\left(\Phi(t,x)-\int_0^1 g_1(z)\Phi(t,(1-z)x)\mathrm{d}z\right)},\quad \theta^{(2)*}(t,x) = e^{-\psi_2\left(\Phi(t,x)-\int_0^1 g_2(z)\Phi(t,z+(1-z)x)\mathrm{d}z\right)}. \tag{20}$$

In this sense, finding $\Phi$ yields the optimal controls. Notice that the interior solutions (20) are obtained with a sufficiently large $\theta_{\max} > 0$ by the uniform bound (11), which we assume in what follows.

## 3. Mathematical analysis and numerical computation

### 3.1 Key Lemmas

This section presents key lemmas on the terms in the right-hand sides of (17). The proof of **Lemma 1** is not presented because it can be checked directly. This lemma provides a basic result important both for the viscosity and distribution solution approaches. It would be important to recall that $\Omega$ is a compact interval, meaning that any $u \in C(\Omega)$ is uniformly continuous in $\Omega$ and belongs to $u \in L^{\infty}(\Omega)$.

***Lemma 1***

*For each* $x \in \Omega$, $H(x,\cdot)$ *is Lipschitz continuous in* $\mathbb{R}$ *with the Lipschitz constant* $L_H a(x)$ with $L_H = r_{\max} + \sigma\lambda_{\max}$.

***Lemma 2***

*For each* $\Phi_1, \Phi_2 \in C(\Omega)$, *there is a constant* $C > 0$ *such that*

$$\left|K(\Phi_1) - K(\Phi_2)\right| \le C\left\|\Phi_1 - \Phi_2\right\|_{C(\Omega)} \quad \textit{in } \Omega. \tag{21}$$

**(Proof of Lemma 2)**

Choose $\Phi_1, \Phi_2 \in C(\Omega)$. We have a constant $L > 0$ such that $|\Phi_1|, |\Phi_2| \le L$ because $\Omega$ is compact. For each $x \in \Omega$, we have

$$\frac{\psi_1}{\nu_1}\left|K(\Phi_1) - K(\Phi_2)\right| = \left| e^{-\psi_1\left(\Phi_1(x)-\int_0^1 g_1(z)\Phi_1((1-z)x)\mathrm{d}z\right)} - e^{-\psi_1\left(\Phi_2(x)-\int_0^1 g_1(z)\Phi_2((1-z)x)\mathrm{d}z\right)} \right|. \tag{22}$$

Then, we get $C = 2\nu_1 e^{2L\psi_1}$ because

$$\frac{\psi_1}{\nu_1}\left|K(\Phi_1)-K(\Phi_2)\right| \le \psi_1 e^{2L\psi_1}\left|\begin{array}{l}\Phi_1(x)-\Phi_2(x)\\ +\left(\int_0^1 g_1(z)\Phi_2((1-z)x)\mathrm{d}z-\int_0^1 g_1(z)\Phi_1((1-z)x)\mathrm{d}z\right)\end{array}\right|$$

$$\begin{aligned}&\le \psi_1 e^{2L\psi_1}\left(\left\|\Phi_1-\Phi_2\right\|_{C(\Omega)}+\int_0^1 g_1(z)\left|\Phi_2((1-z)x)-\Phi_1((1-z)x)\right|\mathrm{d}z\right)\\ &\le 2\psi_1 e^{2L\psi_1}\left\|\Phi_1-\Phi_2\right\|_{C(\Omega)}\end{aligned} . \quad (23)$$

□

***Lemma 3***

*For each* $\Phi_1,\Phi_2 \in C(\Omega)$*, there is a constant* $C>0$ *such that*

$$\left|M(\Phi_1)-M(\Phi_2)\right| \le C\left\|\Phi_1-\Phi_2\right\|_{C(\Omega)},\quad x\in\Omega . \quad (24)$$

The proof of **Lemma 3** is omitted because it is essentially the same with that of **Lemma 2**. **Lemmas 2** and **3** show the mapping property $K,M : C(\Omega)\to C(\Omega)$, with which we can define a (continuous) viscosity solution as presented in **Definition 1**.

The following two lemmas are key for the analysis in the distribution sense.

***Lemma 4***

*For each* $u\in L^\infty(\Omega)$ *, each* $\Phi_1,\Phi_2\in C(\Omega)$ *, any* $\varepsilon>0$ *, there is a constant* $C>0$ *(independent from* $\varepsilon$ *) such that*

$$\left|\int_0^1 u(x)\left(K(\Phi_1)(x)-K(\Phi_2)(x)\right)\mathrm{d}x\right| \le C\int_0^1\left(\varepsilon\left|u(x)\right|^2+\frac{1}{\varepsilon}\left|\Phi_2(x)-\Phi_1(x)\right|^2\right)\mathrm{d}x . \quad (25)$$

**(Proof of Lemma 4)**

A straightforward calculation like that in the proof of **Lemma 2** shows

$$\begin{aligned}&\left|\int_0^1 u(x)\left(K(\Phi_1)(x)-K(\Phi_2)(x)\right)\mathrm{d}x\right|\\ &\le \nu_1 e^{2\psi_1 L}\int_0^1\int_0^1 g_1(z)\left(\left|u(x)\right|\left|\Phi_2(x)-\Phi_1(x)\right|+\left|u(x)\right|\left|\Phi_2((1-z)x)-\Phi_1((1-z)x)\right|\right)\mathrm{d}z\mathrm{d}x\end{aligned} . \quad (26)$$

For any $\varepsilon>0$, the first integral in the right-hand side of (26) is evaluated with the classical Cauchy-Schwartz inequality as

$$\begin{aligned}\int_0^1\int_0^1 g_1(z)\left|u(x)\right|\left|\Phi_2(x)-\Phi_1(x)\right|\mathrm{d}z\mathrm{d}x &= \int_0^1\left|u(x)\right|\left|\Phi_2(x)-\Phi_1(x)\right|\mathrm{d}x\\ &\le \int_0^1\frac{1}{2}\left(\varepsilon\left|u(x)\right|^2+\frac{1}{\varepsilon}\left|\Phi_2(x)-\Phi_1(x)\right|^2\right)\mathrm{d}x\end{aligned} . \quad (27)$$

On the second integral in the right-hand side of (26), we have

$$\begin{aligned}
&\int_0^1\int_0^1 g_1(z)\left|u(x)\right|\left|\Phi_2\left((1-z)x\right)-\Phi_1\left((1-z)x\right)\right|\mathrm{d}z\mathrm{d}x \\
&\le \frac{1}{2}\int_0^1\int_0^1 g_1(z)\left(\varepsilon\left|u(x)\right|^2+\frac{1}{\varepsilon}\left|\Phi_2\left((1-z)x\right)-\Phi_1\left((1-z)x\right)\right|^2\right)\mathrm{d}z\mathrm{d}x \quad . \qquad (28)\\
&= \frac{\varepsilon}{2}\int_0^1\left|u(x)\right|^2\mathrm{d}x+\frac{1}{2\varepsilon}\int_0^1\int_0^1 g_1(z)\left|\Phi_2\left((1-z)x\right)-\Phi_1\left((1-z)x\right)\right|^2\mathrm{d}z\mathrm{d}x
\end{aligned}$$

A change of variables $y=(1-z)x$ gives

$$\begin{aligned}
\int_0^1\int_0^1 g_1(z)\left|\Phi_2\left((1-z)x\right)-\Phi_1\left((1-z)x\right)\right|^2\mathrm{d}z\mathrm{d}x &= \int_0^1 g_1(z)\frac{1}{1-z}\int_0^{1-z}\left|\Phi_2(y)-\Phi_1(y)\right|^2\mathrm{d}y\mathrm{d}z \\
&\le \left(\int_0^1\frac{g(z)}{1-z}\mathrm{d}z\right)\left(\int_0^1\left|\Phi_2(x)-\Phi_1(x)\right|^2\mathrm{d}x\right)
\end{aligned} \quad . \qquad (29)$$

Consequently, we have

$$\begin{aligned}
&\left|\int_0^1 u(x)\left(K(\Phi_1)(x)-K(\Phi_2)(x)\right)\mathrm{d}x\right| \\
&\le \nu_1 e^{2\psi_1 L}\int_0^1\left(\varepsilon\left|u(x)\right|^2+\frac{1}{2\varepsilon}\left(1+\int_0^1\frac{g_1(z)}{1-z}\mathrm{d}z\right)\left|\Phi_2(x)-\Phi_1(x)\right|^2\right)\mathrm{d}x
\end{aligned} \quad . \qquad (30)$$

Choosing $C=\left(1+\int_0^1\frac{g_1(z)}{1-z}\mathrm{d}z\right)\nu_1 e^{2\psi_1 L}$ completes the proof.

□

***Lemma 5***

*For each $u\in L^\infty(\Omega)$, each $\Phi_1,\Phi_2\in C(\Omega)$, any $\varepsilon>0$, there is a constant $C>0$ (independent from $\varepsilon$) such that*

$$\left|\int_0^1 u(x)\left(K(\Phi_1)(x)-K(\Phi_2)(x)\right)\mathrm{d}x\right|\le C\int_0^1\left(\varepsilon\left|u(x)\right|^2+\frac{1}{\varepsilon}\left|\Phi_2(x)-\Phi_1(x)\right|^2\right)\mathrm{d}x. \qquad (31)$$

The proof of **Lemma 5** is omitted because it is essentially the same with that of **Lemma 4**.

### 3.2 Viscosity solution approach

Viscosity solutions are candidates of appropriate weak solutions to problems of degenerate parabolic types, especially for those arising in control problems. Our HJBI equation is not an exception. Its viscosity solutions are defined following Azimzadeh et al. [25] and Øksendal and Sulem [11] where test functions do not appear in the nonlocal terms.

***Definition 1***

*Viscosity sub-solution*

*A function* $\Psi \in C(D) \cap USC(\bar{D})$ *such that* $\Psi(x,T) \le 0$ *for* $x \in \Omega$ *is a viscosity sub-solution if the following conditions are satisfied; for each* $(\bar{t},\bar{x}) \in D$,

$$-\frac{\partial \psi}{\partial t} - \frac{\sigma^2}{2} a(x)^2 \frac{\partial^2 \psi}{\partial x^2} - \left[-\gamma_0 x + \gamma_1 (1-x)\right] \frac{\partial \psi}{\partial x} \le H\left(x, \frac{\partial \psi}{\partial x}\right) - K\Psi - M\Psi \quad \text{at} \quad (\bar{t},\bar{x}) \qquad (32)$$

*for any test functions* $\psi \in C^2(\bar{D})$ *such that* $\Psi - \psi$ *is locally maximized at* $(\bar{t},\bar{x})$ *and* $\Psi \le \psi$ *in a neighborhood of this point.*

*Viscosity super-solution*

*A function* $\Psi \in C(D) \cap LSC(\bar{D})$ *such that* $\Psi(x,T) \ge 0$ *for* $x \in \Omega$ *is a viscosity super-solution if the following conditions are satisfied; for each* $(\bar{t},\bar{x}) \in D$,

$$-\frac{\partial \psi}{\partial t} - \frac{\sigma^2}{2} a(x)^2 \frac{\partial^2 \psi}{\partial x^2} - \left[-\gamma_0 x + \gamma_1 (1-x)\right] \frac{\partial \psi}{\partial x} \ge H\left(x, \frac{\partial \psi}{\partial x}\right) - K\Psi - M\Psi \quad \text{at} \quad (\bar{t},\bar{x}) \qquad (33)$$

*for any test functions* $\phi \in C^2(\bar{D})$ *such that* $\Psi - \psi$ *is locally minimized at* $(\bar{t},\bar{x})$ *and* $\Psi \ge \psi$ *in a neighborhood of this point.*

*Viscosity solution*

*A function* $\Psi \in C(\bar{D})$ *is a viscosity solution if it is a viscosity sub-solution as well as a viscosity super-solution.*

We show uniqueness of viscosity solutions to the HJBI equation (17) through a comparison argument.

***Theorem 1***

*The HJBI equation (17) admits at most one viscosity solution.*

**(Proof of Theorem 1)**

The proof here follows that of Proposition 3.8 of Yoshioka and Yoshioka [34]. The proof below is more involved because the present model has nonlinear and nonlocal terms. The variable transformation $t \to T - t$ is applied to rewrite the terminal value problem as an initial value problem, which fits to the argument in the literature. It is sufficient to show that for any couple of a viscosity sub-solution $\underline{\Phi}$ and a viscosity super-solution $\bar{\Phi}$, we have $\bar{\Phi} \ge \underline{\Phi}$ in $D$. In addition, the transformation $\Phi(t,x) \to \Phi(t,x)e^{Rt}$ with a constant $R > 0$ is used. Then, in the HJBI equation (17), we should formally replace $f(x)$ by $f(x)e^{-Rt}$ and modify other terms accordingly as well. The modified $H, K, M$ (not time-dependent in a continuous manner) are represented in the same notations.

Set $\underline{\Phi}_\delta = \underline{\Phi} - \dfrac{\delta}{T-t}$ with a sufficiently small $\delta > 0$. Clearly, we have $\underline{\Phi}_\delta - \bar{\Phi} \leq 0$ for $t = T$. Assume $\sup_D \left\{ \underline{\Phi}_\delta - \bar{\Phi} \right\} > 0$. Set $f_\varepsilon(s,t,x,y) = \underline{\Phi}_\delta(s,x) - \bar{\Phi}(t,y) - \varphi_\varepsilon(s,x,x,y)$, $\varphi_\varepsilon(s,t,x,y) = \dfrac{1}{2\varepsilon}\left((x-y)^2 + (s-t)^2\right)$ in $\bar{D} \times \bar{D}$. This $f_\varepsilon$ attains a maximum in $D \times D$ because it is upper semi-continuous and $\delta > 0$. The maximizer of $f_\varepsilon$ is denoted as $(s_\varepsilon, t_\varepsilon, x_\varepsilon, y_\varepsilon) \in D \times D$. Then, we have

$$f_\varepsilon(s_\varepsilon, t_\varepsilon, x_\varepsilon, y_\varepsilon) \geq f_\varepsilon(s,s,x,x) = \underline{\Phi}_\delta(s,x) - \bar{\Phi}(s,x) \quad \text{for } (s,x) \in D. \tag{34}$$

Following the standard argument of the doubling of variables [26], we can choose a sequence $\varepsilon = \varepsilon_k$ with $\lim_{k \to +\infty} \varepsilon_k = 0$ such that

$$\lim_{k \to +\infty} s_{\varepsilon_k} = \lim_{k \to +\infty} t_{\varepsilon_k} = s_0, \quad \lim_{k \to +\infty} x_{\varepsilon_k} = \lim_{k \to +\infty} y_{\varepsilon_k} = x_0 \tag{35}$$

$$\lim_{k \to +\infty} \frac{1}{\varepsilon_k}\left(s_{\varepsilon_k} - t_{\varepsilon_k}\right)^2 = 0, \quad \lim_{k \to +\infty} \frac{1}{\varepsilon_k}\left(x_{\varepsilon_k} - y_{\varepsilon_k}\right)^2 = 0 \tag{36}$$

with some $(s_0, x_0) \in D$ such that $\underline{\Phi}_\delta(s_0, x_0) - \bar{\Phi}(s_0, x_0) = \sup_D \left\{ \underline{\Phi}_\delta - \bar{\Phi} \right\} > 0$. Hereafter, we only consider such a sub-sequence with sufficiently large $k$. We have

$$\underline{\Phi}_\delta(s_0, x_0) - \bar{\Phi}(s_0, x_0) \geq \underline{\Phi}_\delta(s,x) - \bar{\Phi}(s,x) \quad \text{for } (s,x) \in D. \tag{37}$$

Rearranging the equation yields

$$\bar{\Phi}(s,x) - \underline{\Phi}_\delta(s,x) \geq \bar{\Phi}(s_0, x_0) - \underline{\Phi}_\delta(s_0, x_0). \tag{38}$$

We see that $\underline{\Phi}_\delta(s,x) - \left(\bar{\Phi}(t_\varepsilon, y_\varepsilon) + \varphi_\varepsilon(s,x,t_\varepsilon,y_\varepsilon)\right)$ is maximized at $(s_\varepsilon, x_\varepsilon)$ and $\bar{\Phi}(t,y) - \left(\underline{\Phi}_\delta(s_\varepsilon, x_\varepsilon) - \varphi_\varepsilon(s_\varepsilon, x_\varepsilon, t, y)\right)$ is minimized at $(t_\varepsilon, y_\varepsilon)$. We then use $\underline{\Phi}_\delta(s_\varepsilon, x_\varepsilon) + \varphi_\varepsilon(s,x,t_\varepsilon,y_\varepsilon) - \varphi_\varepsilon(s_\varepsilon, x_\varepsilon, t_\varepsilon, y_\varepsilon)$ as a test function for the viscosity sub-solution and $\bar{\Phi}(t_\varepsilon, y_\varepsilon) - \varphi_\varepsilon(s_\varepsilon, x_\varepsilon, t, y) + \varphi_\varepsilon(s_\varepsilon, x_\varepsilon, t_\varepsilon, y_\varepsilon)$ as a test function for the viscosity super-solution, respectively.

Set $p_\varepsilon = \dfrac{x_\varepsilon - y_\varepsilon}{\varepsilon}$ and $q_\varepsilon = \dfrac{s_\varepsilon - t_\varepsilon}{\varepsilon}$. By Ishiis' lemma [26], there exist $M_1, M_2$ such that

$$-\frac{3}{\varepsilon}\begin{pmatrix}1 & 0\\ 0 & 1\end{pmatrix} \le \begin{pmatrix}M_1 & 0\\ 0 & -M_2\end{pmatrix} \le \frac{3}{\varepsilon}\begin{pmatrix}1 & -1\\ -1 & 1\end{pmatrix} \tag{39}$$

with

$$\begin{aligned} & q_\varepsilon + R\underline{\Phi}_\delta(s_\varepsilon, x_\varepsilon) - \frac{\sigma^2}{2} a(x_\varepsilon)^2 M_1 - \left[-\gamma_0 x_\varepsilon + \gamma_1(1-x_\varepsilon)\right] p_\varepsilon \\ & \le H(s_\varepsilon, x_\varepsilon, p_\varepsilon) - (K+M)(\underline{\Phi}_\delta)(s_\varepsilon, x_\varepsilon) \end{aligned}, \tag{40}$$

and

$$\begin{aligned} & q_\varepsilon + R\bar{\Phi}(t_\varepsilon, y_\varepsilon) - \frac{\sigma^2}{2} a(y_\varepsilon)^2 M_2 - \left[-\gamma_0 y_\varepsilon + \gamma_1(1-y_\varepsilon)\right] p_\varepsilon \\ & \ge H(t_\varepsilon, y_\varepsilon, p_\varepsilon) - (K+M)(\bar{\Phi})(t_\varepsilon, y_\varepsilon) \end{aligned}. \tag{41}$$

Combining (40) and (41) yields

$$\begin{aligned} & R\left(\underline{\Phi}_\delta(s_\varepsilon, x_\varepsilon) - \bar{\Phi}(t_\varepsilon, y_\varepsilon)\right) + (K+M)(\underline{\Phi}_\delta)(s_\varepsilon, x_\varepsilon) - (K+M)(\bar{\Phi})(t_\varepsilon, y_\varepsilon) \\ & \le H(s_\varepsilon, x_\varepsilon, p_\varepsilon) - H(t_\varepsilon, y_\varepsilon, p_\varepsilon) + \frac{\sigma^2}{2} a(x_\varepsilon)^2 M_1 - \frac{\sigma^2}{2} a(y_\varepsilon)^2 M_2 \end{aligned}. \tag{42}$$

By **Lemma 1** and the structure condition [35], which is fulfilled by the HJBI equation with the transformation of variables, there exists a non-negative and continuous function $c$ in $\Omega$ with $c(0)=0$ such that the right-hand side of (42) is bounded by $c(\varepsilon p_\varepsilon(1+p_\varepsilon)) + O(|t_\varepsilon - s_\varepsilon|)$. Then, we get

$$\begin{aligned} & R\left(\underline{\Phi}_\delta(s_\varepsilon, x_\varepsilon) - \bar{\Phi}(t_\varepsilon, y_\varepsilon)\right) + (K+M)(\underline{\Phi}_\delta)(s_\varepsilon, x_\varepsilon) - (K+M)(\bar{\Phi})(t_\varepsilon, y_\varepsilon) \\ & \le c(\varepsilon p_\varepsilon(1+p_\varepsilon)) + O(|t_\varepsilon - s_\varepsilon|) \end{aligned}. \tag{43}$$

Letting $\varepsilon \to +0$ in (40) and (41) yields $\varepsilon p_\varepsilon(1+p_\varepsilon) \to 0$ by (36), and thus

$$R\left(\underline{\Phi}_\delta - \bar{\Phi}\right) + (K+M)\underline{\Phi}_\delta - (K+M)\bar{\Phi} \le 0 \quad \text{at} \quad (s_0, x_0). \tag{44}$$

Because we have $\bar{\Phi}(s_0, x_0) - \underline{\Phi}_\delta(s_0, x_0) < 0$ and $R > 0$, (44) leads to

$$(K+M)\underline{\Phi}_\delta - (K+M)\bar{\Phi} < 0 \quad \text{at} \quad (s_0, x_0). \tag{45}$$

Now, we prove $(K+M)\underline{\Phi}_\delta \ge (K+M)\bar{\Phi}$ at $(s_0, x_0)$. Firstly, we prove $K\underline{\Phi}_\delta \ge K\bar{\Phi}$ at $(s_0, x_0)$. Assume $K\underline{\Phi}_\delta < K\bar{\Phi}$ at $(s_0, x_0)$. Then, we should have

$$\bar{\Phi}(s_0, x_0) - \int_0^1 g_1(z)\bar{\Phi}(s_0, (1-z)x_0)\mathrm{d}z > \underline{\Phi}_\delta(s_0, x_0) - \int_0^1 g_1(z)\underline{\Phi}_\delta(s_0, (1-z)x_0)\mathrm{d}z\,. \tag{46}$$

We have

$$
\begin{aligned}
&\int_0^1 g_1(z)\bar{\Phi}(s_0,(1-z)x_0)\mathrm{dz}-\int_0^1 g_1(z)\underline{\Phi}_\delta(s_0,(1-z)x_0)\mathrm{dz}\\
&=\int_0^1 g_1(z)\left\{\bar{\Phi}(s_0,(1-z)x_0)-\underline{\Phi}_\delta(s_0,(1-z)x_0)\right\}\mathrm{dz}\\
&\geq\int_0^1 g_1(z)\left\{\bar{\Phi}(s_0,x_0)-\underline{\Phi}_\delta(s_0,x_0)\right\}\mathrm{dz}\\
&=\bar{\Phi}(s_0,x_0)-\underline{\Phi}_\delta(s_0,x_0)
\end{aligned}
\tag{47}
$$

by $\int_0^1 g_1(z)\mathrm{dz}=1$. Substituting (47) into (46) gives the contradiction

$$
\bar{\Phi}(s_0,x_0)-\underline{\Phi}_\delta(s_0,x_0)>\bar{\Phi}(s_0,x_0)-\underline{\Phi}_\delta(s_0,x_0)\,. \tag{48}
$$

Thus, $K\underline{\Phi}_\delta\geq K\bar{\Phi}$ at $(s_0,x_0)$. Similarly, we have $M\underline{\Phi}_\delta\geq M\bar{\Phi}$ at $(s_0,x_0)$. Consequently, we have $(K+M)\underline{\Phi}_\delta\geq(K+M)\bar{\Phi}$ at this point, which contradicts with (45). Therefore, we have $\underline{\Phi}_\delta(s_0,x_0)-\bar{\Phi}(s_0,x_0)\leq 0$. Because $(s_0,x_0)\in D$ is arbitrary, we get the desired result: $\underline{\Phi}_\delta-\bar{\Phi}\leq 0$ in $D$. Finally, we take the limit $\delta\to+0$.

□

**Theorem 1** is useful in convergence of numerical solutions to monotone, stable, and consistent numerical methods [25, 30] toward the viscosity solution by the fundamental theorem of convergence [36]. The following result is a byproduct of **Theorem 1**, from which the uniqueness of the HJBI equation (17) follows.

***Corollary 1***

For any viscosity sub-solution $\underline{\Phi}$ and viscosity super-solution $\bar{\Phi}$, $\bar{\Phi}\geq\underline{\Phi}$ in $D$.

### 3.3 Distribution solution approach

Unique solvability of the HJBI equation (17) from the viewpoint of distribution solutions, which are alternative weak solutions defined in an appropriate Banach space with a weighted norm, is presented. Our result is non-trivial because of using the specific weighted Sobolev space, which is utilized to handle the degenerate diffusion coefficient. In general, standard Sobolev spaces, like $H^1$, are not appropriate for degenerate parabolic problems even in linear cases [28] because of the difficulty to have coercivity of the corresponding bilinear form. In addition, analyzing distribution solutions to the HJBI equation may open the door to its numerical approximation with finite element and related schemes.

The discussion here is based on the approach with the weighted Sobolev space (Proposition 12.1 of Bensoussan [27]) combined with a weight that harmonizes with the degenerate diffusion coefficient (Oleinik and Radkevič [28]).

Set $\mathcal{W} = \{g : g \in C^1(\Omega)\}$. For any $v, w \in \mathcal{W}$, set the scalar product

$$(w,v)_{\mathcal{W}_0} = \int_\Omega \left(a^2(x) w_x v_x + wv\right) \mathrm{d}x + (wv)_{x=0} + (wv)_{x=1} . \tag{49}$$

Define $\mathcal{W}_0$ as a Hilbert space obtained by closing $\mathcal{W}$ with respect to the norm $(w,w)_{\mathcal{W}_0}$. Here, the notation $\frac{\partial w}{\partial x} = w_x$ was utilized. Consider the function space $\mathcal{H} = \left\{\varphi \middle| \varphi \in L^2\left(0,T;\mathcal{W}_0\right), \varphi_t \in L^2\left(0,T;\mathcal{W}_0^*\right)\right\}$ with a free parameter $\mu > 0$ (See, Eq. (50)) and $\mathcal{W}_0^*$ is the dual of $\mathcal{W}_0$, which is chosen later. This space is equipped with the norm $\|\cdot\|$:

$$\|v\|^2 = \int_0^T e^{\mu t} \int_\Omega \left(a^2(x) v_x^2 + v^2\right) \mathrm{d}x\mathrm{d}t + \int_0^T e^{\mu t} \left(\left(v^2\right)_{x=0} + \left(v^2\right)_{x=1}\right) \mathrm{d}t \quad \text{for} \quad v \in \mathcal{H} . \tag{50}$$

Dependence of the norm $\|\cdot\|$ on the parameter $\mu$ is suppressed here for the sake of simplicity of description. We assume $\gamma_0\gamma_1 \neq 0$. The mathematical results in this sub-section can be straightforwardly applied to the case $\gamma_0\gamma_1 = 0$ by omitting the corresponding terms from the equations and the norms following Oleinik and Radkevič [28].

Hereafter, set $\psi_1 = \psi_2 = \psi$ and $\nu_1 = \nu_2 = \nu$ for the sake of simplicity of presentation. The case, $\psi_1 \neq \psi_2$ and/or $\nu_1 \neq \nu_2$ can be handled without technical difficulties. Based on the boundedness result (11), distribution solutions to the HJBI equation (17) are defined as follows.

***Definition 2***

*A bounded function $u \in \mathcal{H}$ is a distribution solution if*

$$\begin{aligned}
&\int_0^T e^{\mu t} \int_\Omega \left(-\frac{\partial u}{\partial t}\right) v \mathrm{d}x\mathrm{d}t + \int_0^T e^{\mu t} \int_\Omega \frac{\sigma^2}{2} a^2(x) u_x v_x \mathrm{d}x\mathrm{d}t \\
&+ \int_0^T e^{\mu t} \int_\Omega \left\{-\left(\gamma_1 - (\gamma_0 + \gamma_1)x\right) \frac{\partial u}{\partial x}\right\} v \mathrm{d}x\mathrm{d}t \\
&= \int_0^T e^{\mu t} \int_\Omega H(x, u_x) v \mathrm{d}x\mathrm{d}t - \int_0^T e^{\mu t} \int_\Omega (Ku) v \mathrm{d}x\mathrm{d}t - \int_0^T e^{\mu t} \int_\Omega (Mu) v \mathrm{d}x\mathrm{d}t
\end{aligned} \quad \textit{for all } v \in \mathcal{H} . \tag{51}$$

The following theorem is our second main result, which gives a second characterization of solutions to the HJBI equation (17). The assumption $a \in C^2(\Omega)$ assumed in the theorem is common in logistic type models [2, 32].

***Theorem 2***

*Assume $a \in C^2(\Omega)$. Then, the HJBI equation (17) subject to the terminal condition $\Phi(T,x) = 0$, $x \in \Omega$ admits at most one distribution solution for sufficiently large $\mu > 0$.*

**(Proof of Theorem 2)**

The proof is based on a contraction mapping argument. Without any loss of generality, set $\sigma^2 = 1$. The calculation here can be justified with a conventional density argument. Firstly, set the map $\Xi : \mathcal{H} \to \mathcal{H}$, $\Xi(v) = u$ for $v \in \mathcal{H}$, as the solution $u \in \mathcal{H}$ to the linear PDIE

$$\begin{aligned} & -\frac{\partial u}{\partial t} - \frac{1}{2}a^2\frac{\partial^2 u}{\partial x^2} - \left(\gamma_1 - (\gamma_0 + \gamma_1)x\right)\frac{\partial u}{\partial x}, \quad (t,x) \in D \\ & = H(x, v_x) - (K+M)(v) \end{aligned} \tag{52}$$

with $u(T,x) = 0$ in $\Omega$. This map $\Xi$ is well-defined. This is because of the unique solvability of the linear equation

$$-\frac{\partial u}{\partial t} - \frac{1}{2}a^2\frac{\partial^2 u}{\partial x^2} - \left(\gamma_1 - (\gamma_0 + \gamma_1)x\right)\frac{\partial u}{\partial x} = f(x), \quad (t,x) \in D, \tag{53}$$

with $u(T,x) = 0$ in $\Omega$ follows from the argument with the Garding's inequality.

To prove the statement, it is sufficient to show that $\Xi$ is a contraction mapping for sufficiently large $\mu > 0$. For $v_1, v_2 \in \mathcal{H}$, set $\Xi(v_1) = u_1$, $\Xi(v_2) = u_2$, $\tilde{v} = v_1 - v_2$, and $\tilde{u} = u_1 - u_2$. We have $\tilde{u}(T,x) = 0$ for $x \in \Omega$ and

$$\begin{aligned} & -\frac{\partial \tilde{u}}{\partial t} - \frac{1}{2}a^2\frac{\partial^2 \tilde{u}}{\partial x^2} - \left(\gamma_1 - (\gamma_0 + \gamma_1)x\right)\frac{\partial \tilde{u}}{\partial x} \\ & = H(x, v_{1,x}) - H(x, v_{2,x}) - (Kv_2 - Kv_1) - (Mv_2 - Mv_1) \end{aligned}, \quad (t,x) \in D. \tag{54}$$

Testing the equation (54) with $\tilde{u}$ yields

$$\begin{aligned} & \int_0^T e^{\mu t}\int_\Omega \left(-\tilde{u}_t - \frac{1}{2}a^2\tilde{u}_{xx} - \left(\gamma_1 - (\gamma_0 + \gamma_1)x\right)\tilde{u}_x\right)\tilde{u}\,\mathrm{d}x\mathrm{d}t \\ & = \int_0^T e^{\mu t}\int_\Omega \left[H(x, v_{1,x}) - H(x, v_{2,x})\right]\tilde{u}\,\mathrm{d}x\mathrm{d}t \\ & -\int_0^T e^{\mu t}\int_\Omega (Kv_2 - Kv_1)\tilde{u}\,\mathrm{d}x\mathrm{d}t - \int_0^T e^{\mu t}\int_\Omega (Mv_2 - Mv_1)\tilde{u}\,\mathrm{d}x\mathrm{d}t \end{aligned}. \tag{55}$$

The left-hand side of (55) is calculated as

$$\begin{aligned} & \int_0^T e^{\mu t}\int_\Omega \left(-\tilde{u}_t - \frac{1}{2}a^2\tilde{u}_{xx} - \left(\gamma_1 - (\gamma_0 + \gamma_1)x\right)\tilde{u}_x\right)\tilde{u}\,\mathrm{d}x\mathrm{d}t \\ & = \int_0^T e^{\mu t}\int_\Omega \left(-\frac{1}{2}\left(\tilde{u}^2\right)_t\right)\mathrm{d}x\mathrm{d}t + \int_0^T e^{\mu t}\int_\Omega \left(-\frac{1}{2}a^2\tilde{u}_{xx}\tilde{u}\right)\mathrm{d}x\mathrm{d}t \\ & +\int_0^T e^{\mu t}\int_\Omega \left(-\left(\gamma_1 - (\gamma_0 + \gamma_1)x\right)\tilde{u}_x\tilde{u}\right)\mathrm{d}x\mathrm{d}t \end{aligned}. \tag{56}$$

The first term in (56) becomes

$$\int_0^T e^{\mu t}\int_\Omega\left(-\frac{1}{2}\left(\tilde{u}^2\right)_t\right)\mathrm{d}x\mathrm{d}t=\int_0^T\int_\Omega\left(-\frac{1}{2}\left(\tilde{u}^2\right)_t\right)e^{\mu t}\mathrm{d}x\mathrm{d}t$$
$$=\int_\Omega\left(-\frac{1}{2}\left(\tilde{u}^2\right)_{t=T}e^{\mu T}+\frac{1}{2}\left(\tilde{u}^2\right)_{t=0}\right)\mathrm{d}x+\mu\int_0^T\int_\Omega\frac{1}{2}\tilde{u}^2e^{\mu t}\mathrm{d}x\mathrm{d}t\,. \tag{57}$$
$$=\frac{1}{2}\int_\Omega\left(\tilde{u}^2\right)_{t=0}\mathrm{d}x+\frac{\mu}{2}\int_0^T e^{\mu t}\int_\Omega\tilde{u}^2\mathrm{d}x\mathrm{d}t$$

By an integration by parts, the second term in (56) becomes

$$\int_0^T e^{\mu t}\int_\Omega\left(-\frac{1}{2}a^2\tilde{u}_{xx}\tilde{u}\right)\mathrm{d}x\mathrm{d}t=\int_0^T e^{\mu t}\int_\Omega\left(\frac{1}{2}a^2\tilde{u}\right)_x\tilde{u}_x\mathrm{d}x\mathrm{d}t$$
$$=\frac{1}{2}\int_0^T e^{\mu t}\int_\Omega a^2\tilde{u}_x^2\mathrm{d}x\mathrm{d}t+\int_0^T e^{\mu t}\int_\Omega aa_x\tilde{u}\tilde{u}_x\mathrm{d}x\mathrm{d}t$$
$$=\frac{1}{2}\int_0^T e^{\mu t}\int_\Omega a^2\tilde{u}_x^2\mathrm{d}x\mathrm{d}t+\int_0^T e^{\mu t}\int_\Omega aa_x\left(\frac{1}{2}\tilde{u}^2\right)_x\mathrm{d}x\mathrm{d}t \qquad . \tag{58}$$
$$=\frac{1}{2}\int_0^T e^{\mu t}\int_\Omega a^2\tilde{u}_x^2\mathrm{d}x\mathrm{d}t-\frac{1}{2}\int_0^T e^{\mu t}\int_\Omega\left(a_xa_x+aa_{xx}\right)\tilde{u}^2\mathrm{d}x\mathrm{d}t$$

The third term in (56) becomes

$$\int_0^T e^{\mu t}\int_\Omega\left(-\left(\gamma_1-\left(\gamma_0+\gamma_1\right)x\right)\tilde{u}_x\right)\tilde{u}\mathrm{d}x\mathrm{d}t$$
$$=\int_0^T e^{\mu t}\int_\Omega\left(-\left(\gamma_1-\left(\gamma_0+\gamma_1\right)x\right)\frac{1}{2}\left(\tilde{u}^2\right)_x\right)\mathrm{d}x\mathrm{d}t \qquad . \tag{59}$$
$$=\frac{1}{2}\int_0^T e^{\mu t}\left[\gamma_0\left(\tilde{u}^2\right)_{x=1}+\gamma_1\left(\tilde{u}^2\right)_{x=0}\right]\mathrm{d}t-\frac{1}{2}\int_0^T e^{\mu t}\int_\Omega\left(\gamma_0+\gamma_1\right)\tilde{u}^2\mathrm{d}x\mathrm{d}t$$

By **Lemmas 4** and **5**, there exists a global constant $L>0$ such that

$$\int_0^T e^{\mu t}\int_\Omega\left(Kv_2-Kv_1\right)\tilde{u}\mathrm{d}x\mathrm{d}t+\int_0^T e^{\mu t}\int_\Omega\left(Mv_2-Mv_1\right)\tilde{u}\mathrm{d}x\mathrm{d}t\le L\int_0^T e^{\mu t}\int_\Omega\left(\varepsilon\tilde{u}^2+\frac{1}{\varepsilon}\tilde{v}^2\right)\mathrm{d}x\mathrm{d}t \tag{60}$$

with an arbitrary constant $\varepsilon>0$. In addition, by **Lemma 1**, we have

$$\int_0^T e^{\mu t}\int_\Omega\left[H\left(x,v_{1,x}\right)-H\left(x,v_{2,x}\right)\right]\left|\tilde{u}\right|\mathrm{d}x\mathrm{d}t\le L_H\int_0^T e^{\mu t}\int_\Omega a\left|v_{1,x}-v_{2,x}\right|\left|\tilde{u}\right|\mathrm{d}x\mathrm{d}t$$
$$=L_H\int_0^T e^{\mu t}\int_\Omega a\left|\tilde{v}_x\right|\left|\tilde{u}\right|\mathrm{d}x\mathrm{d}t \qquad . \tag{61}$$
$$\le\frac{L_H}{2}\int_0^T e^{\mu t}\int_\Omega\left(\frac{1}{\varepsilon}a^2\left|\tilde{v}_x\right|^2+\varepsilon\tilde{u}^2\right)\mathrm{d}x\mathrm{d}t$$

Therefore, we get

$$\frac{1}{2}\int_\Omega\left(\tilde{u}^2\right)_{t=0}\mathrm{d}x+\frac{1}{2}\int_0^T e^{\mu t}\int_\Omega a^2\tilde{u}_x^2\mathrm{d}x\mathrm{d}t+\int_0^T e^{\mu t}\int_\Omega\frac{1}{2}\left(\mu-a_xa_x-aa_{xx}-\left(\gamma_0+\gamma_1\right)\right)\tilde{u}^2\mathrm{d}x\mathrm{d}t$$
$$+\frac{1}{2}\int_0^T e^{\mu t}\left[\gamma_0\left(\tilde{u}^2\right)_{x=1}+\gamma_1\left(\tilde{u}^2\right)_{x=0}\right]\mathrm{d}t \tag{62}$$
$$\le L\int_0^T e^{\mu t}\int_\Omega\left(\varepsilon\tilde{u}^2+\frac{1}{\varepsilon}\tilde{v}^2\right)\mathrm{d}x\mathrm{d}t+\frac{L_H}{2}\int_0^T e^{\mu t}\int_\Omega\left(\frac{1}{\varepsilon}a^2\left|\tilde{v}_x\right|^2+\varepsilon\tilde{u}^2\right)\mathrm{d}x\mathrm{d}t$$

and thus

$$\begin{aligned}&\frac{1}{2}\int_0^T e^{\mu t}\int_\Omega a^2\tilde{u}_x^2 \mathrm{d}x\mathrm{d}t+\frac{1}{2}\int_0^T e^{\mu t}\left[\gamma_0\left(\tilde{u}^2\right)_{x=1}+\gamma_1\left(\tilde{u}^2\right)_{x=0}\right]\mathrm{d}t\\&+\frac{1}{2}\int_0^T e^{\mu t}\int_\Omega\left(\mu-2L\varepsilon-L_H\varepsilon-a_xa_x-aa_{xx}-\left(\gamma_0+\gamma_1\right)\right)\tilde{u}^2\mathrm{d}x\mathrm{d}t\,.\\&\le\frac{L}{\varepsilon}\int_0^T e^{\mu t}\int_\Omega\tilde{v}^2\mathrm{d}x\mathrm{d}t+\frac{L_H}{2}\frac{1}{\varepsilon}\int_0^T e^{\mu t}\int_\Omega a^2\left|\tilde{v}_x\right|^2\mathrm{d}x\mathrm{d}t\end{aligned}\tag{63}$$

The right-hand side of (63) is estimated from above as

$$\frac{L}{\varepsilon}\int_0^T e^{\mu t}\int_\Omega\tilde{v}^2\mathrm{d}x\mathrm{d}t+\frac{L_H}{2}\frac{1}{\varepsilon}\int_0^T e^{\mu t}\int_\Omega a^2\left|\tilde{v}_x\right|^2\mathrm{d}x\mathrm{d}t\le\frac{L'}{\varepsilon}\left\|\tilde{v}\right\|^2\tag{64}$$

with $L'=\max\{L,L_H/2\}$. On the other hand, the left-hand side of (63) is bounded from below as

$$\begin{aligned}&\frac{1}{2}\int_0^T e^{\mu t}\int_\Omega a^2\tilde{u}_x^2 \mathrm{d}x\mathrm{d}t+\frac{1}{2}\int_0^T e^{\mu t}\int_\Omega\left(\mu-2L\varepsilon-L_H\varepsilon-a_xa_x-aa_{xx}-\left(\gamma_0+\gamma_1\right)\right)\tilde{u}^2\mathrm{d}x\mathrm{d}t\\&+\frac{1}{2}\int_0^T e^{\mu t}\left[\gamma_0\left(\tilde{u}^2\right)_{x=1}+\gamma_1\left(\tilde{u}^2\right)_{x=0}\right]\mathrm{d}t\\&\ge\frac{1}{2}\int_0^T e^{\mu t}\int_\Omega a^2\tilde{u}_x^2 \mathrm{d}x\mathrm{d}t+\frac{1}{2}\left(\mu-C_\varepsilon\right)\int_0^T e^{\mu t}\int_\Omega\tilde{u}^2\mathrm{d}x\mathrm{d}t\\&+\frac{1}{2}\int_0^T e^{\mu t}\left[\gamma_0\left(\tilde{u}^2\right)_{x=1}+\gamma_1\left(\tilde{u}^2\right)_{x=0}\right]\mathrm{d}t\end{aligned}\tag{65}$$

with $C_\varepsilon=2L\varepsilon+L_H\varepsilon+\max_\Omega\left(a_xa_x+aa_{xx}\right)+\left(\gamma_0+\gamma_1\right)$. Set $\varepsilon=4L'\min\{1,\gamma_0,\gamma_1\}^{-1}$ and choose a sufficiently large $\mu>0$ such that $\mu-C_\varepsilon\ge 1$. Consequently, we have

$$\begin{aligned}&\frac{1}{2}\int_0^T e^{\mu t}\int_\Omega a^2\tilde{u}_x^2 \mathrm{d}x\mathrm{d}t+\frac{1}{2}\int_0^T e^{\mu t}\int_\Omega\left(\mu-2L\varepsilon-L_H\varepsilon-a_xa_x-aa_{xx}-\left(\gamma_0+\gamma_1\right)\right)\tilde{u}^2\mathrm{d}x\mathrm{d}t\\&+\frac{1}{2}\int_0^T e^{\mu t}\left[\gamma_0\left(\tilde{u}^2\right)_{x=1}+\gamma_1\left(\tilde{u}^2\right)_{x=0}\right]\mathrm{d}t\\&\ge\frac{1}{2}\int_0^T e^{\mu t}\int_\Omega a^2\tilde{u}_x^2 \mathrm{d}x\mathrm{d}t+\frac{1}{2}\int_0^T e^{\mu t}\int_\Omega\tilde{u}^2\mathrm{d}x\mathrm{d}t+\frac{1}{2}\int_0^T e^{\mu t}\left[\gamma_0\left(\tilde{u}^2\right)_{x=1}+\gamma_1\left(\tilde{u}^2\right)_{x=0}\right]\mathrm{d}t\\&\ge\frac{1}{2}\min\{1,\gamma_0,\gamma_1\}\left\|\tilde{u}\right\|^2\end{aligned}\tag{66}$$

and thus

$$\left\|\tilde{u}\right\|^2\le\frac{1}{2}\left\|\tilde{v}\right\|^2.\tag{67}$$

The inequality (67) shows $\Xi$ is a contraction in $\mathscr{H}$ with the above-specified $\mu>0$.

□

The following corollary is the byproduct of **Theorem 2** showing boundedness of distribution solutions.

***Corollary 2***

*Under the assumption of* ***Theorem 2****, for sufficiently large* $\mu$*, we have the explicit bound*

$$\|u\| \le C f_{\max} . \tag{68}$$

**(Proof of Corollary 2)**

As in the proof of **Theorem 2**, we have

$$\begin{aligned}
&\int_0^T e^{\mu t}\int_\Omega \left(-u_t - \frac{1}{2}a^2 u_{xx} - \left(\gamma_1 - \left(\gamma_0+\gamma_1\right)x\right)u_x\right) u \mathrm{d}x\mathrm{d}t \\
&= \frac{1}{2}\int_\Omega \left(u^2\right)_{t=0}\mathrm{d}x + \frac{1}{2}\int_0^T e^{\mu t}\left[\gamma_0\left(u^2\right)_{x=1} + \gamma_1\left(u^2\right)_{x=0}\right]\mathrm{d}t \\
&+\frac{1}{2}\int_0^T e^{\mu t}\int_\Omega a^2u_x^2\mathrm{d}x\mathrm{d}t + \frac{1}{2}\int_0^T e^{\mu t}\int_\Omega\left[\mu - \left(a_x a_x + a a_{xx}\right) - \left(\gamma_0+\gamma_1\right)\right]u^2\mathrm{d}x\mathrm{d}t
\end{aligned} . \tag{69}$$

By **Lemmas 4** and **5**, there exists a constant $\bar{L} > 0$ such that

$$\begin{aligned}
&\int_0^T e^{\mu t}\int_\Omega H\left(x,u_x\right)u\mathrm{d}x\mathrm{d}t - \int_0^T e^{\mu t}\int_\Omega\left(\left(K+M\right)u\right)u\mathrm{d}x\mathrm{d}t \\
&\le \frac{\bar{L}}{2}\int_0^T e^{\mu t}\int_\Omega\left(\varepsilon u^2 + \frac{1}{\varepsilon}a^2u_x^2\right)\mathrm{d}x\mathrm{d}t + \bar{L}\int_0^T e^{\mu t}\int_\Omega u^2\mathrm{d}x\mathrm{d}t + f_{\max}\|u\|
\end{aligned} \tag{70}$$

with arbitrary $\varepsilon > \bar{L}$. Combining (69) and (70) leads to

$$\begin{aligned}
&\frac{1}{2}\int_\Omega\left(\tilde{u}^2\right)_{t=0}\mathrm{d}x + \frac{1}{2}\int_0^T e^{\mu t}\left[\gamma_0\left(u^2\right)_{x=1}+\gamma_1\left(u^2\right)_{x=0}\right]\mathrm{d}t + \frac{1}{2}\left(1-\frac{\bar{L}}{\varepsilon}\right)\int_0^T e^{\mu t}\int_\Omega a^2u_x^2\mathrm{d}x\mathrm{d}t \\
&+\frac{1}{2}\int_0^T e^{\mu t}\int_\Omega\left[\mu - \left(a_xa_x + aa_{xx}\right) - \left(\gamma_0+\gamma_1\right) - \bar{L}\left(\varepsilon+2\right)\right]u^2\mathrm{d}x\mathrm{d}t \\
&\le f_{\max}\|u\|
\end{aligned} . \tag{71}$$

By choosing $\varepsilon = 2\bar{L}$, the left-hand side of (71) is bounded from below as

$$\begin{aligned}
&\frac{1}{2}\int_\Omega\left(\tilde{u}^2\right)_{t=0}\mathrm{d}x + \frac{1}{2}\int_0^T e^{\mu t}\left[\gamma_0\left(u^2\right)_{x=1}+\gamma_1\left(u^2\right)_{x=0}\right]\mathrm{d}t + \frac{1}{4}\int_0^T e^{\mu t}\int_\Omega a^2u_x^2\mathrm{d}x\mathrm{d}t \\
&+\frac{1}{2}\int_0^T e^{\mu t}\int_\Omega\left[\mu - \left(a_xa_x+aa_{xx}\right) - \left(\gamma_0+\gamma_1\right) - 2\bar{L}\left(\bar{L}+1\right)\right]u^2\mathrm{d}x\mathrm{d}t \\
&\ge \frac{1}{2}\min\left\{\gamma_0,\gamma_1,\frac{1}{2}\right\}\left(\int_0^T e^{\mu t}\int_\Omega b^2u_x^2\mathrm{d}x\mathrm{d}t + \int_0^T e^{\mu t}\left[\left(u^2\right)_{x=1}+\left(u^2\right)_{x=0}\right]\mathrm{d}t\right) \\
&+\frac{1}{2}\int_0^T e^{\mu t}\int_\Omega\left[\mu - \left(a_xa_x+aa_{xx}\right) - \left(\gamma_0+\gamma_1\right) - 2\bar{L}\left(\bar{L}+1\right)\right]u^2\mathrm{d}x\mathrm{d}t \\
&\ge \frac{1}{2}\min\left\{\gamma_0,\gamma_1,\frac{1}{2}\right\}\|u\|^2
\end{aligned} \tag{72}$$

with sufficiently large $\mu$ such that

$$\mu \ge \max_\Omega\left(a_xa_x + aa_{xx}\right) + \left(\gamma_0+\gamma_1\right) + 2\bar{L}\left(\bar{L}+1\right) + \min\left\{\gamma_0,\gamma_1,\frac{1}{2}\right\} . \tag{73}$$

With such a $\mu$, combining (71) with (72) yields

$$\|u\| \le 2 f_{\max}\min\left\{\gamma_0,\gamma_1,\frac{1}{2}\right\}^{-1} . \tag{74}$$

Choosing $C = 2 f_{\max}\min\left\{\gamma_0,\gamma_1,\frac{1}{2}\right\}^{-1}$ completes the proof.

□

Another important outcome of **Theorem 2** is the following comparison principle, characterizing dependence of the solution to the parameters and coefficients. Assume that $H = H(x,p)$ continuously depends on a real parameter $\omega \in \mathbb{R}$ for each $(x,p) \in \Omega \times \mathbb{R}$, and this parameter dependence is expressed with a subscript as $H_\omega$. Assume the monotone dependence $H_{\omega_2} \geq H_{\omega_1}$ for $\omega_1 \leq \omega_2$. Then, we can state the following corollary. The proof is omitted here since it just follows the Proof of Proposition 12.3 of Bensoussan [27].

***Corollary 3***

*The solution to the HJBI equation (17) with $H = H_{\omega_j}$ is denoted as $u_j$ ( $j = 1,2$ , $\omega_1 \leq \omega_2$). Then, $u_2 \geq u_1$.*

The specific form of $H$ further leads to the next result, clearly showing that a less costly and/or less ambiguity-averse performance index leads to s smaller value function.

***Corollary 4***

*The solution to the HJBI equation (17) with is increasing with respect to the parameter $\psi_0$ and functions $f$ and $h$.*

### 3.4 Numerical computation

Demonstrative computational examples of the HJBI equation (17) are presented focusing on a single-species population management problem where the population should be kept near a prescribed value through human interventions [31]. Here, an ergodic limit (well-mixed limit, $T \to +\infty$), which is a topic that has not been mathematically analyzed above, is numerically explored. A particular focus is put on dependence of the value function and optimal controls on the model parameters to complement the mathematical analysis results.

#### 3.4.1 Numerical method

A monotone, stable, and bounded finite difference scheme is used to discretization of the HJBI equation (17). The scheme used here is the maximum use of central difference [37] combined with the linear interpolation discretization of the nonlocal terms [25]. The temporal integration is carried out in a fully-implicit way. **Theorem 1** shows that a necessary (but not sufficient) condition of convergence of numerical solutions toward the viscosity solution locally and uniformly in $D$. To guarantee full convergence, one should check monotonicity, stability, and consistency of the scheme [36]. The discretization methods are monotone and stable, and the scheme is consistent. Then, numerical solutions converge to the value function if the value function is the unique viscosity solution. Formal convergence rate of monotone schemes, including ours, is at most first-order in space-time because the scheme uses a monotone discretization in space and a fully-implicit discretization in time.

The discretized system is solved at each time step with a standard policy iteration algorithm [30]. Each iteration in the algorithm solves one linear system based on the Thomas method [38]. The domain $\Omega$ is discretized into 500 cells with 501 vertices. The temporal integration of the HJBI equation (17) is carried out backward in time with the time increment $0.005$. This computational resolution has been found to be sufficiently fine for the analysis below. The computation is carried out from the terminal time $t = T$ toward the initial time $t = 0$ since (17) is a time-backward problem. $\theta_{\max}$ and $\lambda_{\max}$ are set as 100, which are sufficiently large values such that the computed $\lambda^*$, $\theta^{(1)*}$, and $\theta^{(2)*}$ are valued inside of their admissible ranges.

### 3.4.2 Problem without the control $q$

The present model is firstly applied to a simplified problem without the control $q$ chosen by the human ($q \equiv 0$), in which the optimal controls are symmetric with respect to the center $x = 0.5$ of the domain $\Omega$. This problem is considered as a worst-case estimation of the average net disutility. We specify the coefficients as follows: $\sigma = 1$, $a = x(1-x)$, $\nu_1 = \nu_2 = \nu = 1$, $\psi_0 = 0.5$, $\psi_1 = \psi_2 = \psi = 0.5$, $f(x) = \max\{2x-1, 1-2x\}$, $T = 50$, unless otherwise specified. In the present case, the function $f$ implies that the population should be kept near its minimizing point $x = 0.5$ and we have $f_{\max} = 1$. The terminal time $T$ is chosen so that the optimal controls are sufficiently close to a time-independent state (an ergodic state) at the initial time $t = 0$. In such a case, based on the ergodic ansatz [34, 39], the control variables are time-independent and $-\frac{\partial \Phi}{\partial t} \to E = const > 0$ over $\Omega$ as $t$ decreases, where $E$ is the effective Hamiltonian identified as a limit of the time-averaged value function in $\Omega$. The numerical approach taken here is the large $T$ method that has been effectively applied to an environmental control problem [34]. The problem with large $T$ is well-posed by the mathematical analysis results.

**Figure 1** shows the computed normalized value function $\Phi - \min_{\Omega} \Phi$ ($\min_{\Omega} \Phi = 38.665$ in this case) and the optimal controls $\lambda^*$, $\theta^{(1)*}$, and $\theta^{(2)*}$ at $t=0$. The value function $\Phi$ is symmetric with respect to $x=0.5$, and the controls $\theta^{(1)*}$ and $\theta^{(2)*}$ are symmetric with each other. The control $\lambda^*$ has a point symmetry at $x=0.5$. The computed $-\frac{\partial \Phi}{\partial t}$ with a backward difference is sufficiently converged to a constant $E=0.7943$ in $\Omega$ at $t=0$, supporting the ergodic ansatz.

Parameter dependence of the optimal controls and the effective Hamiltonian are explored numerically. **Figures 2** and **3** show $\lambda^*$ for different values of $\psi_0$ and $\theta^{(1)*}$ for different values of $\psi$, respectively, in $\Omega$ at $t=0$. The amplitudes of the controls $\lambda^*$ and $\theta^{(1)*}$ become larger as the potential ambiguity, namely $\psi_0$ and $\psi$, increases as expected. **Figure 4** shows the effective Hamiltonian $E$ for different values of $\psi_0$ and $\psi$. Theoretically, we should have $0 \le E \le f_{\max} = 1$, which is clearly satisfied in the numerical solution. The effective Hamiltonian $E$ is increasing with respect to $\psi_0$ and $\psi$, and continuously approaches toward $f_{\max} = 1$ from below. **Corollary 4** supports the increasing nature of $E$ because it can be seen as a minimum value of the time-averaged $\Phi$.

### 3.4.3 Problem with the control $q$

The present model with the coefficients specified in the previous sub-section is applied to the problem with the control $q$ in which the distributions of the controls are asymmetric. The coefficients $r(q)=1-q$ and $h(q)=0.1\times q$, $q_{\max}=1$ are specified that in turn lead to a bang-bang type optimal control representing the intervention to suppress the growth ($q^*=1$) and no intervention ($q^*=0$). We expect that, at each $t<T$, there exists a sub-domain $\Omega_1$ of $\Omega$ such that $q^*(t,x)=1$ in $\Omega_1$ and $q^*(t,x)=0$ outside of $\Omega_1$, because the interventions decrease the population and should not be carried out at least if the population is small. An interest here is dependence of this sub-domain on the parameters in the performance index $\phi$.

**Figure 5** shows the computed normalized value function $\Phi - \min_{\Omega} \Phi$ ($\min_{\Omega} \Phi = 37.003$ in this case) and the optimal controls $\lambda^*$, $\theta^{(1)*}$, $\theta^{(2)*}$, and $q^*$ at $t=0$. The symmetry of the value function and the optimal controls are now broken due to the non-zero $q^*$. The computational results imply that including the human interventions, despite that it is costly, reduces the value function. The effective Hamiltonian. The computed $-\frac{\partial \Phi}{\partial t}$ is sufficiently converged to a constant $E$ in $\Omega$ at $t=0$ in the present case as well, and $E$ for different values of $\psi_0$ and $\psi$ are included in **Figure 4**. The computational results again comply with the theoretical results as in the problem without $q$. The computed optimal control $q^*$ suggests a bang-bang control having two switching points, with which the population can be effectively confined near the minimum point of $f$. The control should not be activated if the population is sufficiently large or small, suggesting an intensive management policy focusing on the situation when the population is close to the desired state. **Figures 6** and **7** show the dependence of the sub-domain $\Omega_1$ (where $q^*(0,x)=1$) on $\psi_0$ and $\psi$, respectively. The computational results suggest that doing no intervention is optimal if the decision-maker puts a sufficiently small trust on the model. We see that the optimal control $q^*$ is affected not only by the ambiguity of diffusion but also those by jumps, suggesting the non-negligible interactions between the controls and ambiguities. Overall, the computational results in this paper suggest that the model problem has a continuous dependence on the parameters $\psi_0$ and $\psi$ on ambiguity.

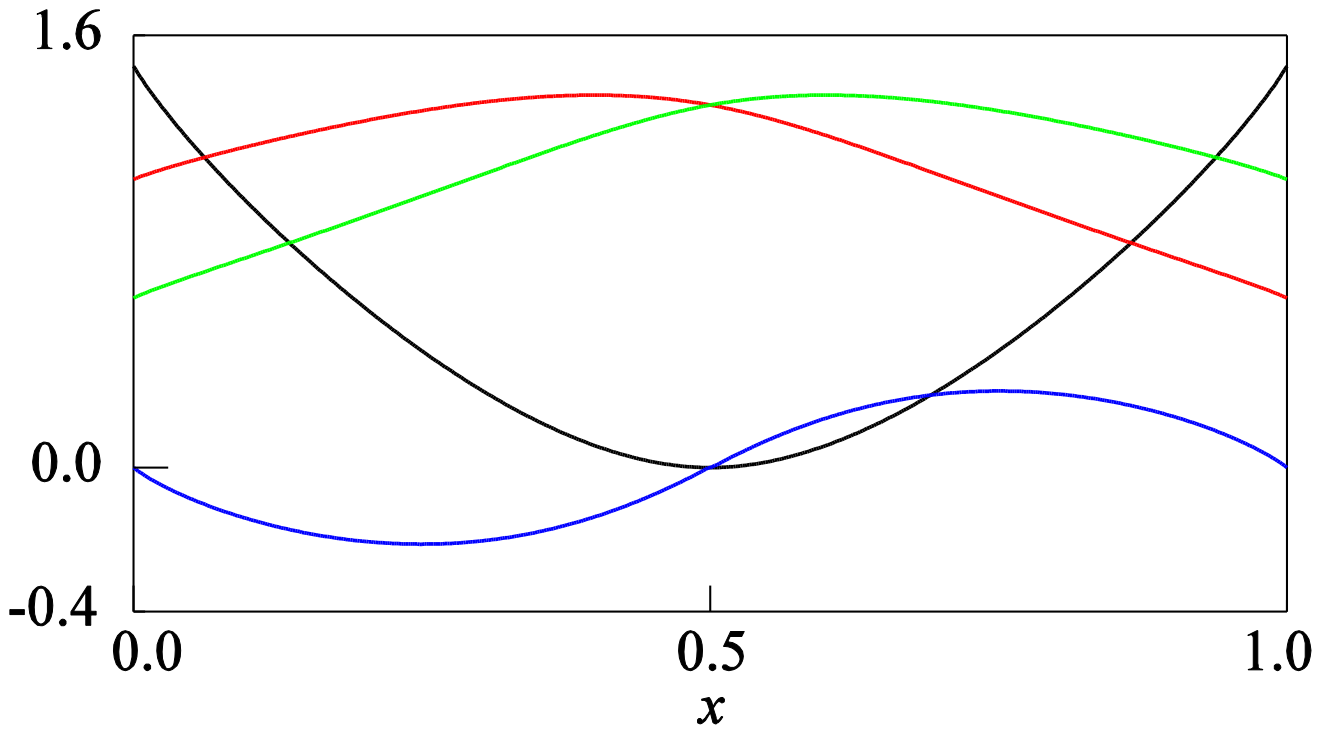

**Figure 1.** The computed $\Phi - \min_{\Omega} \Phi$ (black), $\lambda^*$ (blue), $\theta^{(1)*}$ (red), and $\theta^{(2)*}$ (green).

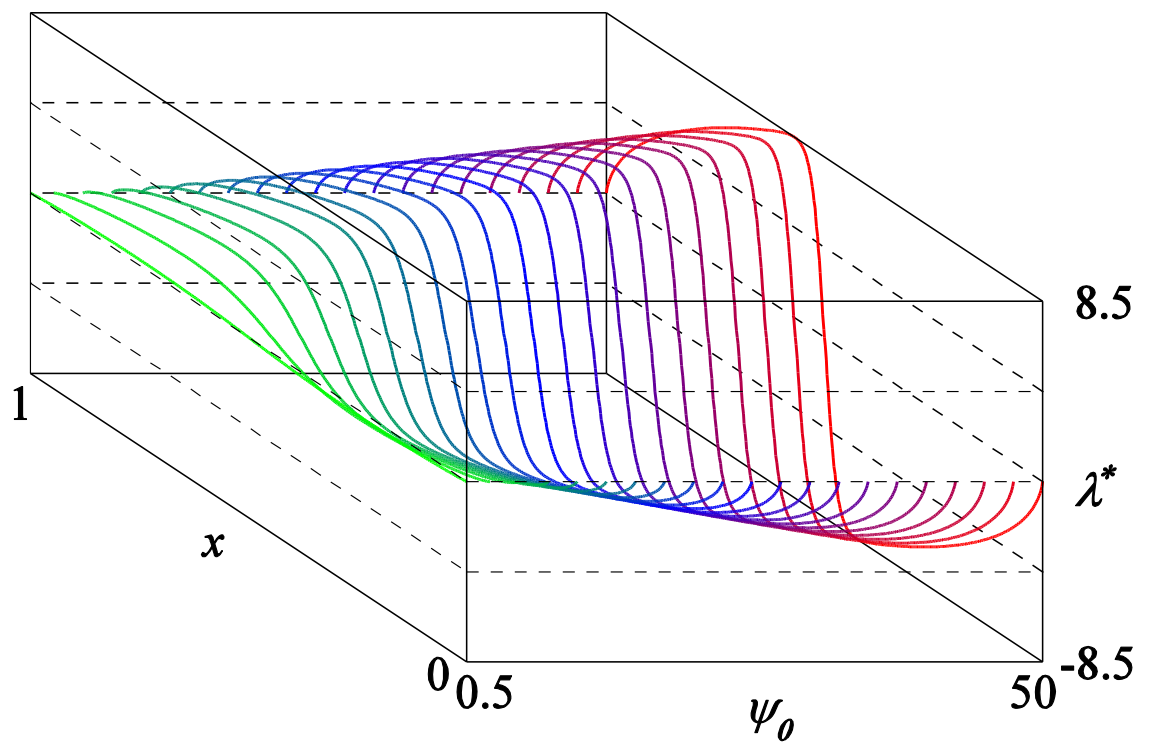

**Figure 2.** Computed $\lambda^*$ for different values of $\psi_0$.

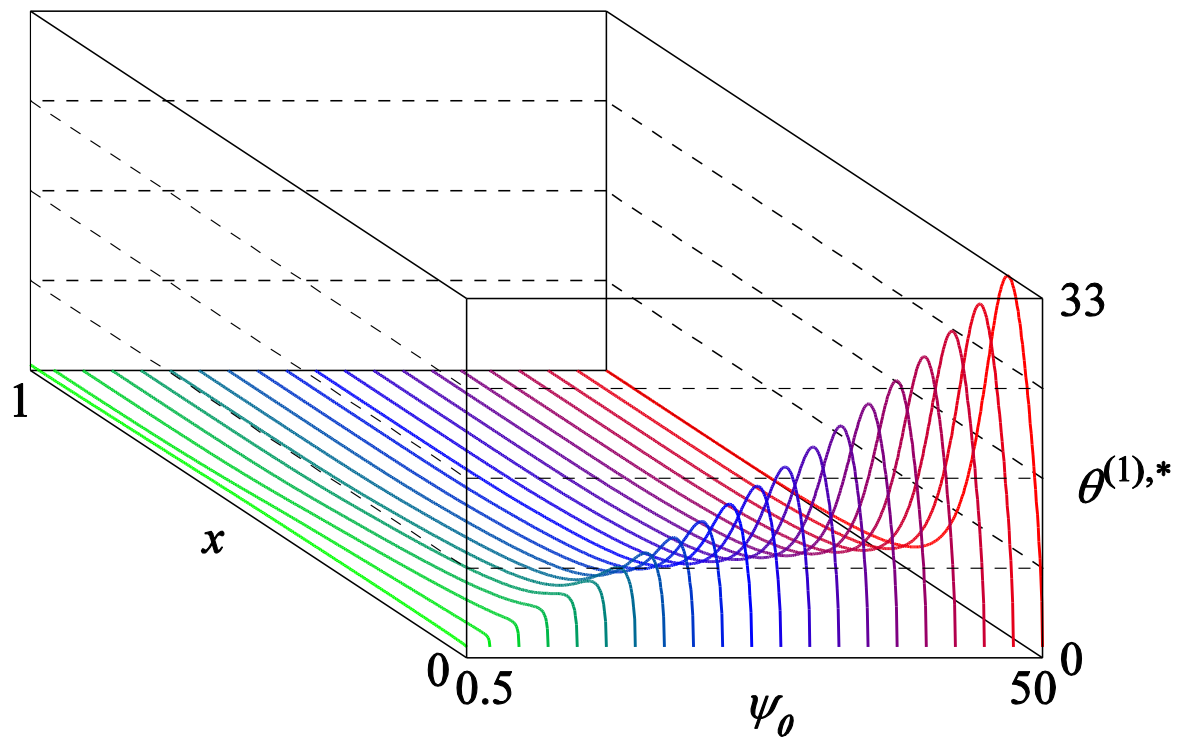

**Figure 3.** Computed $\theta^{(1)*}$ for different values of $\psi_0$.

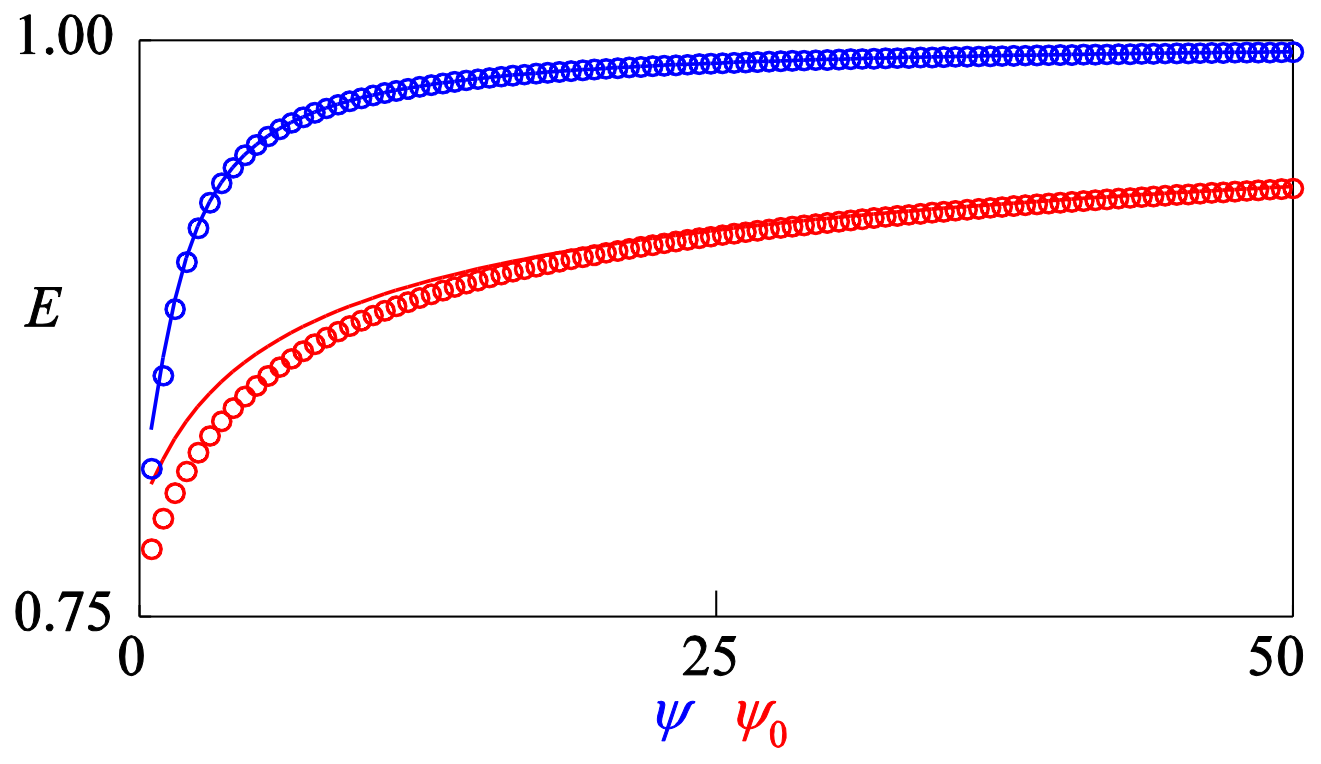

**Figure 4.** Computed $E$ for different values of $\psi_0$ and $\psi$. The solid lines are $E$ for the problem without considering $q$ and circles are that for the problem with $q$.

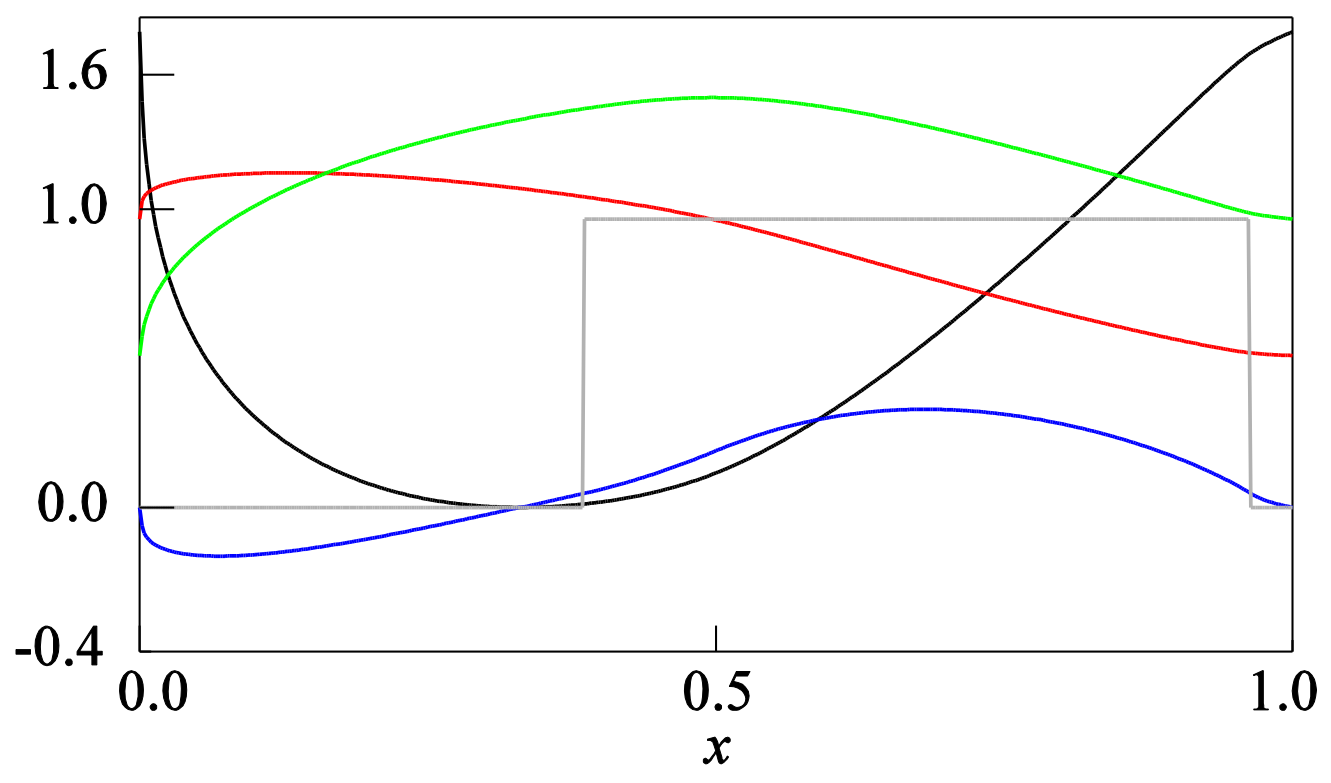

**Figure 5.** The computed $\Phi - \min_{\Omega} \Phi$ (black), $\lambda^*$ (blue), $\theta^{(1)*}$ (red), $\theta^{(2)*}$ (green), and $q^*$ (grey).

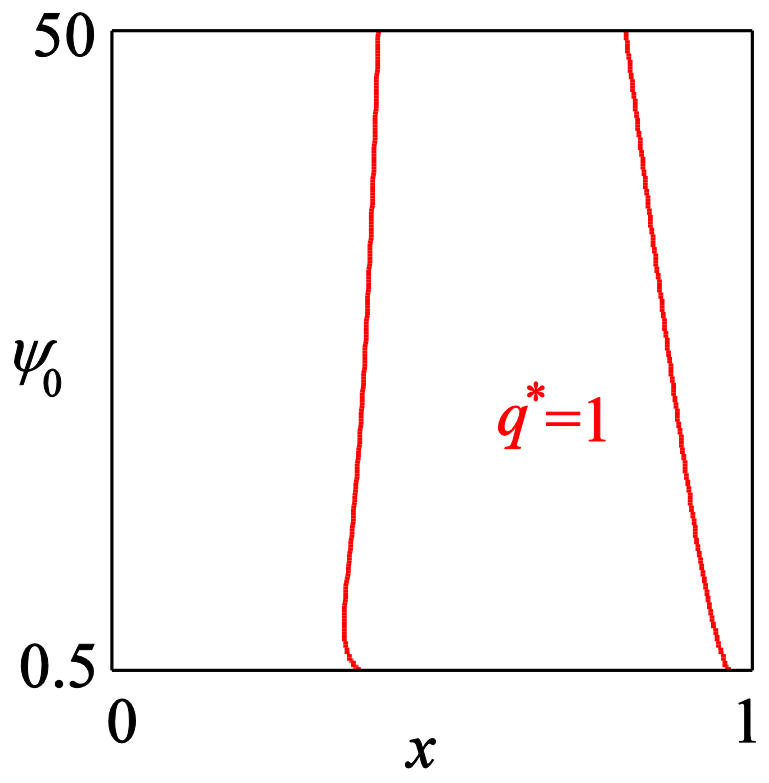

**Figure 6.** Dependence of the sub-domain $\Omega_1$ on $\psi_0$.

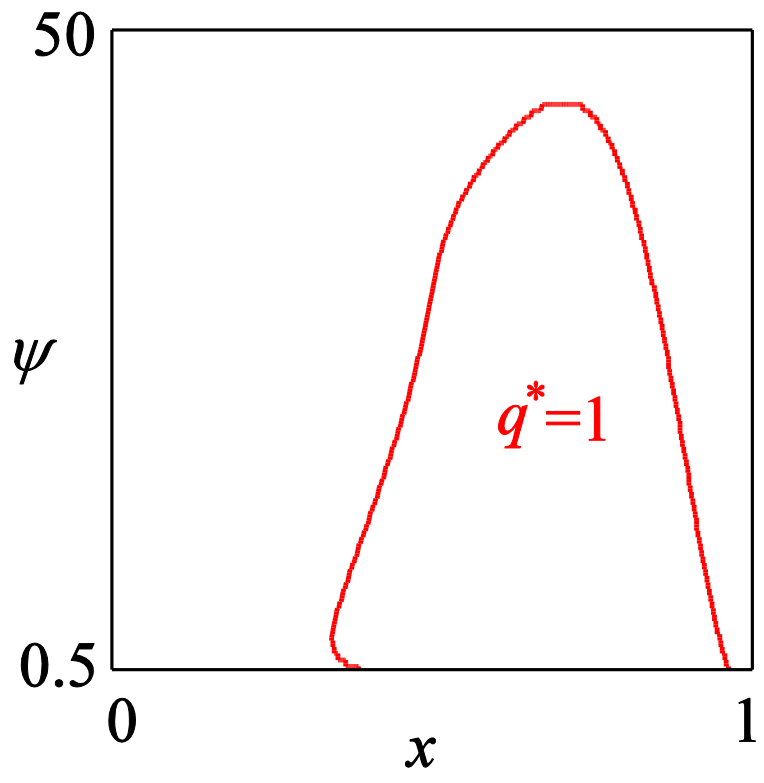

**Figure 7.** Dependence of the sub-domain $\Omega_1$ on $\psi$ .

## 4. Conclusions

An HJBI equation associated with a biological management problem of a bounded and ambiguous stochastic population dynamics was mathematically analyzed from the standpoints of viscosity and distribution solutions. The presented formalism provides an efficient way for optimization of stochastic systems subject to model ambiguity. The problem turned out to be well-posed from both standpoints. It was demonstrated that its numerical solutions can be reasonably computed with a simple finite difference scheme under an ergodic limit.

The presented mathematical analysis can be extended to problems of multi-dimensional jump-diffusion processes with ambiguity where some of the variables are bounded. Mathematical analysis of problems with unbounded ambiguity, which would be subject to a worse regularity of the Hamiltonians, is an open issue to be addressed in future research. Models with more complicated problems where the diffusion part is also controlled [34] are encountered in modern environmental management problems. A more generic bounded jump noise [40, 41] and/or a more realistic drift [42] can be considered in the present modeling framework. Considering the incomplete information structure [5] is interesting as well. Application of specific related models to modeling and control of other biological phenomena, like algae bloom, is undergoing [43]. Their analysis will be a key future research topic. Mean field modeling considering a number of decision-makers (stakeholders) [44] is also an interesting topic. In applications, not only the optimal controls, but also the controlled state dynamics are of importance as well [45], which would require a well-balanced numerical method having satisfactory accuracy. This research topic is currently undergoing. We have assumed a strong dynamic programming principle for our problem, but establishment of a dynamic programming principle is in general not a trivial issue [46, 47]. Verifying this assumption is a remaining issue, especially if the control variables are allowed to be unbounded. Addressing these issues would be able to advance modeling and analysis on environmental and ecological management.

**Acknowledgements**

JSPS Research Grant No. 18K01714 and a grant for ecological survey of a life history of the landlocked Ayu *Plecoglossus altivelis altivelis* from MLIT of Japan support this research. The comments from an anonymous referee lead to improvement of this paper.

**References**

[1] Thieme, H. R. (2018). Mathematics in Population Biology. Princeton University Press, Princeton.

[2] Lungu, E. M., & Øksendal, B. (1997). Optimal harvesting from a population in a stochastic crowded environment. Mathematical Biosciences, 145(1), 47-75. https://doi.org/10.1016/S0025-5564(97)00029-1

[3] Xu, C. (2017). Global threshold dynamics of a stochastic differential equation SIS model. Journal of Mathematical Analysis and Applications, 447(2), 736-757. https://doi.org/10.1016/j.jmaa.2016.10.041

[4] Yang, B., Cai, Y., Wang, K., & Wang, W. (2019). Optimal harvesting policy of logistic population model in a randomly fluctuating environment. Physica A: Statistical Mechanics and its Applications, 526, Paper ID: 120817. https://doi.org/10.1016/j.physa.2019.04.053

[5] Grandits, P., Kovacevic, R. M., & Veliov, V. M. (2019). Optimal control and the value of information for a stochastic epidemiological SIS-model. Journal of Mathematical Analysis and Applications, 476(2), 665-695. https://doi.org/10.1016/j.jmaa.2019.04.005

[6] Yoshioka, H. (2019a). A simplified stochastic optimization model for logistic dynamics with control-dependent carrying capacity. Journal of Biological Dynamics, 13(1), 148-176. https://doi.org/10.1080/17513758.2019.1576927

[7] Yoshioka, H., Yaegashi, Y., Yoshioka, Y., & Hamagami, K. (2019a). Hamilton–Jacobi–Bellman quasi-variational inequality arising in an environmental problem and its numerical discretization. Computers & Mathematics with Applications, 77(8), 2182-2206. https://doi.org/10.1016/j.camwa.2018.12.004

[8] Yoshioka, H., Yaegashi, Y., Yoshioka, Y., & Tsugihashi, K. (2019b). Optimal harvesting

policy of an inland fishery resource under incomplete information. Applied Stochastic Models in Business and Industry, 35(4), 939-962. https://doi.org/10.1002/asmb.2428

[9] d'Onofrio, A. (2013). Bounded noises in physics, biology, and engineering. Springer New York. https://doi.org/10.1007/978-1-4614-7385-5

[10] He, S., Tang, S., & Wang, W. (2019). A stochastic SIS model driven by random diffusion of air pollutants. Physica A: Statistical Mechanics and its Applications, 121759. https://doi.org/10.1016/j.physa.2019.121759

[11] Øksendal, B., & Sulem, A. (2019). Applied Stochastic Control of Jump Diffusions. Springer, Cham. https://doi.org/10.1007/978-3-030-02781-0

[12] Befekadu, G. K., & Zhu, Q. (2018). Optimal control of diffusion processes pertaining to an opioid epidemic dynamical model with random perturbations. Journal of Mathematical Biology, 78, 1425-1438. https://doi.org/10.1007/s00285-018-1314-y

[13] Bonnemain, T., & Ullmo, D. (2019). Mean field games in the weak noise limit: A WKB approach to the Fokker–Planck equation. Physica A: Statistical Mechanics and its Applications, 523, 310-325. https://doi.org/10.1016/j.physa.2019.01.143

[14] Li, H., & Guo, G. (2019). Dynamic decision of transboundary basin pollution under emission permits and pollution abatement. Physica A: Statistical Mechanics and its Applications, 121869. https://doi.org/10.1016/j.physa.2019.121869

[15] Tabibian, B., Upadhyay, U., De, A., Zarezade, A., Schölkopf, B., & Gomez-Rodriguez, M. (2019). Enhancing human learning via spaced repetition optimization. Proceedings of the National Academy of Sciences, 116(10), 3988-3993. https://doi.org/10.1073/pnas.1815156116

[16] Hansen, L., & Sargent, T. J. (2001). Robust control and model uncertainty. American Economic Review, 91(2), 60-66. https://doi.org/10.1257/aer.91.2.60

[17] Ma, J., & Niu, Y. (2019). The timing and intensity of investment under ambiguity. The North American Journal of Economics and Finance, 419, 318-330. https://doi.org/10.1016/j.najef.2019.04.015

[18] Cartea, A., Jaimungal, S., & Qin, Z. (2016). Model uncertainty in commodity markets. SIAM Journal on Financial Mathematics, 7(1), 1-33. https://doi.org/10.1137/15M1027243

[19] Cartea, Á., Jaimungal, S., & Qin, Z. (2018). Speculative trading of electricity contracts in interconnected locations. Energy Economics, 79, 3-20. https://doi.org/10.1016/j.eneco.2018.11.019

[20] Wang, P., & Li, Z. (2018). Robust optimal investment strategy for an AAM of DC pension plans with stochastic interest rate and stochastic volatility. Insurance: Mathematics and

Economics, 80, 67-83. https://doi.org/10.1016/j.insmatheco.2018.03.003

[21] Zhang, Q., & Chen, P. (2019). Robust optimal proportional reinsurance and investment strategy for an insurer with defaultable risks and jumps. Journal of Computational and Applied Mathematics, 356, 46-66. https://doi.org/10.1016/j.cam.2019.01.034

[22] Manoussi, V., Xepapadeas, A., & Emmerling, J. (2018). Climate engineering under deep uncertainty. Journal of Economic Dynamics and Control, 94, 207-224. https://doi.org/10.1016/j.jedc.2018.06.003

[23] Yoshioka H. and Tsujimura M.: A model problem of stochastic optimal control subject to ambiguous jump intensity, The 23rd Annual International Real Options Conference London, UK, June 27-29, 2019. Full paper (29 pp.): http://www.realoptions.org/openconf2019/data/papers/370.pdf.

[24] Yoshioka, H. (2019). A stochastic differential game approach toward animal migration. Theory in Biosciences, 138(2), 277-303. https://doi.org/10.1007/s12064-019-00292-4

[25] Azimzadeh, P., Bayraktar, E., & Labahn, G. (2018). Convergence of Implicit Schemes for Hamilton-Jacobi-Bellman Quasi-Variational Inequalities. SIAM Journal on Control and Optimization, 56(6), 3994-4016. https://doi.org/10.1137/18M1171965

[26] Crandall, M. G., Ishii, H., & Lions, P. L. (1992). User’s guide to viscosity solutions of second order partial differential equations. Bulletin of the American mathematical society, 27(1), 1-67. https://doi.org/10.1090/S0273-0979-1992-00266-5

[27] Bensoussan, A. (2018). Estimation and control of dynamical systems (Vol. 48). New York: Springer. https://doi.org/10.1007/978-3-319-75456-7

[28] Oleinik, O. & Radkevič, E. V. (1973). Second-order equations with nonnegative characteristic form. Springer, Boston. https://doi.org/10.1007/978-1-4684-8965-1

[29] Ishii, H. (1995). On the equivalence of two notions of weak solutions, viscosity solutions and distribution solutions. Funkcial. Ekvac, 38(1), 101-120.

[30] Neilan, M., Salgado, A. J., & Zhang, W. (2017). Numerical analysis of strongly nonlinear PDEs. Acta Numerica, 26, 137-303. https://doi.org/10.1017/S0962492917000071

[31] Yaegashi, Y., Yoshioka, H., Unami, K., & Fujihara, M. (2018). A singular stochastic control model for sustainable population management of the fish-eating waterfowl *Phalacrocorax carbo*. Journal of Environmental Management, 219, 18-27. https://doi.org/10.1016/j.jenvman.2018.04.099

[32] Lv, J., Wang, K., & Jiao, J. (2015). Stability of stochastic Richards growth model. Applied Mathematical Modelling, 39(16), 4821-4827. https://doi.org/10.1016/j.apm.2015.04.016

[33] Hu, D., & Wang, H. (2019). Reinsurance contract design when the insurer is ambiguity-

averse. Insurance: Mathematics and Economics, 86, 241-255. https://doi.org/10.1016/j.insmatheco.2019.03.007

[34] Yoshioka H. and Yoshioka Y.: Modeling stochastic operation of reservoir under ambiguity with an emphasis on river management. Optimal Control Applications and Methods, 40, 764-790. 10.1002/oca.2510.

[35] S. Koike, A Beginner’s Guide to the Theory of Viscosity Solutions, Mathematical Society of Japan, Tokyo (2004).

[36] Barles, G., & Souganidis, P. E. (1991). Convergence of approximation schemes for fully nonlinear second order equations. Asymptotic Analysis, 4(3), 271-283. https://doi.org/10.3233/ASY-1991-4305

[37] Wang, J., & Forsyth, P. A. (2008). Maximal use of central differencing for Hamilton–Jacobi–Bellman PDEs in finance. SIAM Journal on Numerical Analysis, 46(3), 1580-1601. https://doi.org/10.1137/060675186

[38] L. H. Thomas, Elliptic problems in linear difference equations over a network, Watson Sci. Comput. Lab. Rept., Columbia Univ. New York, vol. 1, 1949

[39] Barles, G. & Souganidis, P. E. (2000). On the large time behavior of solutions of Hamilton-Jacobi equations. SIAM J. Math. Anal., 31(4), 925-939. https://doi.org/10.1137/S0036141099350869

[40] Denisov, S. I., & Bystrik, Y. S. (2019). Statistics of bounded processes driven by Poisson white noise. Physica A: Statistical Mechanics and its Applications, 515, 38-46. https://doi.org/10.1016/j.physa.2018.09.158

[41] Denisov, S. I., & Bystrik, Y. S. (2019). Exact stationary solutions of the Kolmogorov-Feller equation in a bounded domain. Communications in Nonlinear Science and Numerical Simulation, 74, 248-259. https://doi.org/10.1016/j.cnsns.2019.03.023

[42] Wang, Y. (2019). Analysis of a budworm growth model with jump-diffusion. Physica A: Statistical Mechanics and its Applications, 121763. https://doi.org/10.1016/j.physa.2019.121763

[43] Yoshioka H., & Tsujimura M. Stochastic control of single-species population dynamics model subject to jump ambiguity. (under review)

[44] Lasry, J. M., & Lions, P. L. (2018). Mean-field games with a major player. Comptes Rendus Mathematique, 356, 886-890. https://doi.org/10.1016/j.crma.2018.06.001

[45] Yaegashi, Y., Yoshioka, H., Tsugihashi, K., & Fujihara, M. (2019). Analysis and computation of probability density functions for a 1-D impulsively controlled diffusion process. Comptes Rendus Mathematique, 357, 306-315.

https://doi.org/10.1016/j.crma.2019.02.007

[46] Koike, S., & Święch, A. (2013). Representation formulas for solutions of Isaacs integro-PDE. Indiana University Mathematics Journal, 1473-1502. https://doi.org/10.1512/iumj.2013.62.5109

[47] Lv, S. (2020). Two-player zero-sum stochastic differential games with regime switching. Automatica, 114, 108819. https://doi.org/10.1016/j.automatica.2020.108819